\documentclass[11pt]{article}

\usepackage[utf8]{inputenc}

\usepackage[
backend=biber,
style=numeric-comp,
sorting=nty
]{biblatex}
\usepackage[margin=1in]{geometry}
\usepackage{csquotes}
\usepackage{amsfonts,amsthm,amsmath,amssymb,accents,mathtools,physics} 
\usepackage{mathrsfs,amsbsy} 
\usepackage{enumitem}
\usepackage{graphicx}
\usepackage{float}
\usepackage{xcolor}
\usepackage{hyperref}
\usepackage{textcomp} 
\usepackage{relsize}

\usepackage{subcaption}
\usepackage{lipsum}
\usepackage{epstopdf}
\usepackage{algorithmic}
\usepackage{forest}
\ifpdf
  \DeclareGraphicsExtensions{.eps,.pdf,.png,.jpg}
\else
  \DeclareGraphicsExtensions{.eps}
\fi

\usepackage{bbm}
\usepackage{dsfont}

\numberwithin{equation}{section}

\usepackage{yhmath}
 \usepackage{booktabs} 

\newtheorem{thm}{Theorem}[section]

\newtheorem{rem}{Remark}[section]
\newtheorem{prop}{Property}[section]

\theoremstyle{definition}

\usepackage{amsopn}

\usepackage{color}

\hypersetup{
    colorlinks=true,
    linkcolor=blue,
    citecolor=blue,
    filecolor=blue,      
    urlcolor=blue,
}

\title{A High-Order Rank-Adaptive Implicit Algorithm for Solving High Dimensional Diffusion Equations using the Hierarchical Tucker Decomposition}

\author{
  Paolo Bosques-Paulet\footnote{This work was done while PBP was an undergraduate student at Swarthmore College. Email: pbosque1@swarthmore.edu. ORCID iD: 0009-0007-3464-7901}
  \and
  Project Advisor: Joseph Nakao\footnote{Department of Mathematics and Statistics, Swarthmore College, Swarthmore, PA, USA. Email: jnakao1@swarthmore.edu. ORCID iD: 0009-0008-6589-4013}
}

\date{}

\begin{document}

\maketitle


\begin{abstract}
\noindent This paper presents a high-order rank-adaptive implicit integrator for the tensor solution of high-dimensional diffusion equations. We extend the 3D version of this method from the Tucker decomposition to higher dimensions using the hierarchical Tucker (HT) decomposition, since the storage complexity for the Tucker decomposition increases exponentially with the number of dimensions $d>3$. The HT format avoids this issue by decomposing the solution according to a binary tree consisting of bases for each dimension and core tensors which connect the bases. Spectral methods are considered for spatial discretization, and diagonally implicit Runge-Kutta methods are considered for time discretization. At each stage of the Runge-Kutta method, the bases computed at the previous stages are augmented to predict the upcoming basis and construct projection subspaces. By projecting onto these enriched subspaces, the bases and cores can be updated in a sequential manner going up the tree from leaf-to-root. Unlike the 3D Tucker method which has a single core tensor, the HT method also updates the intermediate core tensors. Numerical experiments demonstrate that the method observes high-order accuracy, and test how well the integrator captures the solution rank for various sets of time-dependent diffusion coefficients.
\end{abstract}



\section{Introduction}
We develop a high-order implicit rank-adaptive algorithm for solving high-dimensional diffusion equations using the hierarchical Tucker tensor decomposition. In particular, we consider diffusion equations of the form
\begin{equation}
\frac{\partial u}{\partial t}
=
\sum_{i=1}^{d}
D_i(t)\frac{\partial^2u}{\partial x_i^2},
\qquad
\mathbf{x}\in\Omega\subset\mathbb{R}^{d},
\end{equation}
where the diffusion coefficients $D_i(t)$ may vary in time. High-dimensional diffusion equations arise in a wide range of scientific and engineering applications, including stochastic processes, uncertainty quantification, PDE constrained optimization, and kinetic theory, where the evolution of multivariate probability distributions is governed by diffusive dynamics. One of the major computational bottlenecks is the curse of dimensionality. Discretizing the numerical solution with $N$ degrees of freedom in each spatial dimension requires $N^{d}$ unknowns, causing both computational cost and memory requirements to grow exponentially with the number of dimensions. Consequently, conventional grid-based methods become computationally intractable in high dimensions, even for problems of moderate dimensionality. To address this challenge, low-rank tensor representations have emerged as an effective means of exploiting the inherent low-rank structure present in many high-dimensional solutions. Just as the singular value decomposition (SVD) is used to compress 2D matrix solutions, low-rank tensor decompositions can be used to truncate high-dimensional solutions. By approximating multidimensional arrays using compressed tensor formats, low-rank methods significantly reduce storage requirements and computational complexity while retaining the essential features of the solution \cite{Kolda2009,Grasedyck2013}.

High dimensional PDEs  ($d\geq 3$) have been solved using low-rank tensor decompositions including the Tucker, hierarchical Tucker, and tensor train formats; see the review papers \cite{bachmayr2023low,EINKEMMER2025114191}.
A common framework for evolving time-dependent low-rank approximations is the dynamical low-rank (DLR) approach. DLR methods evolve the tensor factors (bases and coefficients) of the solution by projecting onto a corresponding low-rank manifold. The foundational work of Koch and Lubich introduced DLR approximation for matrix-valued problems and established the projection-based formulation for evolving low-rank representations \cite{Koch2007}. This was extended to tensor approximations in \cite{Koch2010}. However, the original formulation is known to be ill-conditioned in some situations, and regularization techniques were used to first address this issue \cite{Kieri2016,Nonnenmacher2008}. The projector-splitting DLR integrator \cite{Lubich2014} and the basis-update and Galerkin (BUG) formulation \cite{BUG_Parallel,Ceruti2022RankAdaptive,Ceruti2022Unconventional} were also developed to improve robustness. Subsequent developments extended these ideas to Tucker tensors and tree tensor network representations \cite{Ceruti2023TTN,Lubich2018,Ceruti2022RankAdaptive,Ceruti2021TTN}.

Implicit time integration is often desired for evolving diffusive dynamics since explicit methods are subject to stability constraints that become increasingly restrictive under mesh refinement. Considerable effort has been devoted to developing implicit rank-adaptive integrators. Much of the work in this direction has been done for the matrix case \cite{Appelo2025Robust,Li2026HighOrder,Rodgers2023Implicit,BUG_Parallel,2d_rail,Naderi2025,Sutti2024,Meng2025,Ceruti2024BUG}, and progress has been made in the higher-dimensional tensor case \cite{Ceruti2021TTN,Lubich2018Tucker,ElKahza,3d_rail,wang2026implicit,Sands2025,ghahremani2026implicit}. Although these approaches have significantly expanded the range of tractable high-dimensional problems that can be solved, there is still a great need for \textit{high-order implicit rank-adaptive methods for dimensions $d\geq 4$}. The work in \cite{Sands2025} solves high-dimensional kinetic simulations using high-order implicit-explicit time-stepping using the hierarchical Tucker format. The recent work in \cite{ghahremani2026implicit} is capable of solving high-dimensional linear and nonlinear equations with high-order BDF methods in the tensor train-cross (TT-cross) format. Yet, given the limited number of high-order implicit rank-adaptive methods to high-dimensional tree tensor networks (such as hierarchical Tucker and tensor train decompositions), there remains a great need for additional developments.

In this work, we extend the Reduced Augmentation Implicit Low-rank (RAIL) framework of \cite{2d_rail,3d_rail} to the $d\geq 4$ case through the hierarchical Tucker decomposition, following a similar idea to \cite{Ceruti2021TTN,Ceruti2023TTN}. The RAIL method falls into the same category of projection-based methods that includes DLR and BUG methods. In three dimensions, the RAIL method stores the spatially discretized solution at each time step in a Tucker decomposition. Each time step updates the three 1D-bases and subsequently evolves the single core tensor storing the coefficients. When integrating using a high-order implicit Runge-Kutta (RK) method, the tensor equation is projected at each stage onto a subspace generated by the bases computed at the previous stages. This process is designed to preserve the high-order temporal accuracy of the RK method. Unlike Tucker-based methods, whose storage requirements continue to grow with dimension $d\geq 4$, the hierarchical Tucker representation uses a dimension tree network to achieve improved scalability. The Tucker and hierarchical Tucker decompositions have comparable storage requirements for $d=3$. Our formulation updates not only the 1D-bases but also the multiple core tensors associated with each intermediate node of the dimension tree during every Runge-Kutta stage. These updated core tensors are subsequently incorporated into the Galerkin projections defining each reduced system. Numerical experiments demonstrate the high-order temporal accuracy, while effectively capturing the evolution of hierarchical Tucker tensor ranks.

The remainder of this paper is organized as follows. Section~2 reviews the hierarchical Tucker decomposition and introduces the notation used throughout the paper. Section~3 develops the hierarchical Tucker reduced augmentation implicit low-rank algorithm. Finally, Section~4 presents numerical experiments demonstrating the accuracy of the proposed method and its ability to capture the evolution of hierarchical Tucker tensor ranks for different sets of diffusion coefficients. Conclusions and ongoing work is described in Section~5.

\section{Preliminaries}
The following are a set of preliminary definitions and properties. They are presented in summary from \cite{KoMa14,Kormann2017LowRankTD}.

\subsection{Tensor Notation and Definitions}
High-dimensional tensors naturally arise when storing the solution of a high-dimensional PDE. For example, the scalar solution $u(x_1,x_2,x_3,x_4)$ is discretized over uniform computational grids in each dimension
\begin{equation} \label{eq:compgrid}
    x_{i,1}<   x_{i,2}<...  <x_{i,{N_i}},\qquad i\in\{1,2,3,4\}.
\end{equation}

We store the values of $u(x_1,x_2,x_3,x_4)$ evaluated over the Cartesian tensor product of the 1D grids in the 4D array (order-4 tensor) $\mathcal{U}\in \mathbb{R}^{N_1\times N_2\times N_3 \times N_4}$. Since we do not assume the reader has experience in low-rank tensors, we first review key terminology and tensor properties. In this paper, we will use the notation from \cite{Kormann2017LowRankTD,KoMa14}. Tensors, or n-dimensional arrays, will be denoted with uppercase caligraphic letters $\mathcal{A}$.
Matrices will be notated in bold-face capital letters $\mathbf{A}$. The subscript of a matrix or tensor will refer to the dimensions it contains information about. However, for the identity matrix $\mathbf{I}$ the subscript will be taken to denote its size.

\subsection*{Mode-$n$ fibers}
The mode-$n$ fibers of a tensor fix all but the $n$th index. The continuous analogue is fixing a function in all but a single variable.

\subsection*{Mode-$n$ matricization}
We can arrange the mode-$n$ fibers of a tensor $\mathcal{A}$ as the columns of a matrix. This is known as the mode-$n$ matricization/flattening/unfolding of the tensor $\mathcal{A}$, denoted by $\mathbf{A}_{(n)}\equiv\text{mat}_n(\mathcal{A})$. We note that the same idea extends to flattening in more then one dimension.

\subsection*{Vectorization}
The vectorization of a tensor, denoted $\operatorname{vec}(\mathcal{A})\in\mathbb{R}^{PQr}$, is the vector formed by arranging all the mode-$n$ fibers into a single vector, for a chosen $n$. The ordering of the rearranged columns in the vectorization does not matter as long as it is consistent.

\subsection*{Kronecker Product for Matrices}
The Kronecker product of two matrices $\mathbf{A} \in \mathbb{R}^{I \times J}$ and 
$\mathbf{B} \in \mathbb{R}^{P \times Q}$, denoted by $\mathbf{A} \otimes \mathbf{B} \in \mathbb{R}^{IP \times JQ}$, 
is the matrix defined by

\begin{equation}
\mathbf{A} \otimes \mathbf{B} =
\begin{pmatrix}
a_{11}\mathbf{B} & \cdots & a_{1J}\mathbf{B} \\
\vdots  & \ddots & \vdots  \\
a_{I1}\mathbf{B} & \cdots & a_{IJ}\mathbf{B}
\end{pmatrix}
=
\big[
a_{:,1} \otimes b_{:,1},\;
a_{:,1} \otimes b_{:,2},\;
\ldots,\;
a_{:,J} \otimes b_{:,Q}
\big].
\end{equation}

\begin{prop}\label{prop1} $
\mathbf{AB} \otimes \mathbf{CD} = (\mathbf{A} \otimes \mathbf{C})(\mathbf{B} \otimes \mathbf{D}).
$
\end{prop}
\begin{prop}\label{prop2}
$\operatorname{vec}(\mathbf{A}\mathbf{X}\mathbf{B}^{T}) = (\mathbf{B}\otimes \mathbf{A})\operatorname{vec}(\mathbf{X}).$
\end{prop}

\subsection*{Tensor mode-$n$ product}

For simplicity, consider a third order tensor $\mathcal{A} \in \mathbb{R}^{P \times Q \times r}$ and matrices 
$\mathbf X \in \mathbb{R}^{P' \times P}$, 
$\mathbf Y \in \mathbb{R}^{Q' \times Q}$, and 
$\mathbf Z \in \mathbb{R}^{r' \times r}$. Although standard notation uses numbers to indicate the mode \cite{Kolda2009}, we use lower-case letters to avoid confusion later in the paper. The mode-$a$ product between the tensor $\mathcal{A}$ and the matrix $\mathbf X$, denoted 
$\mathcal{A} \times_a \mathbf X$, is a tensor of size 
$P' \times Q \times r$. Similarly, the mode-$b$ product is $\mathcal{A} \times_b \mathbf Y$, and the mode-$c$ product is $\mathcal{A} \times_c \mathbf Z$. Elementwise,
\begin{align}
(\mathcal{A} \times_a \mathbf X)_{i'jk}
=
\sum_{i=1}^{P} a_{ijk}\, x_{i'i},
\qquad
i' = 1,\dots,P', \;
j = 1,\dots,Q, \;
k = 1,\dots,r.\\
(\mathcal{A} \times_b \mathbf Y)_{ij'k}
=
\sum_{j=1}^{Q} a_{ijk}\, y_{j'j},
\qquad
i = 1,\dots,P, \;
j' = 1,\dots,Q', \;
k = 1,\dots,r. \\
(\mathcal{A} \times_c \mathbf Z)_{ijk'}
=
\sum_{k=1}^{r} a_{ijk}\, z_{k'k},
\qquad
i = 1,\dots,P, \;
j = 1,\dots,Q, \;
k' = 1,\dots,r'.
\end{align}

Equivalently stated, the mode-$a$ product between $\mathcal{A}$ and $\mathbf{X}$ multiplies 
each mode-$a$ fiber of $\mathcal{A}$ by the matrix $\mathbf{X}$. In the matricized form,
\begin{align}
\mathcal{H} = \mathcal{A} \times_a \mathbf  X
  &\Longleftrightarrow  
\mathbf{H}_{(a)} =\mathbf X \mathbf A_{(a)}, 
\\
\mathcal{J} = \mathcal{A} \times_b \mathbf Y
   &\Longleftrightarrow  
\mathbf{J}_{(b)} =\mathbf  Y\mathbf  A_{(b)},  \\
\mathcal{K} = \mathcal{A} \times_c \mathbf  Z
  &\Longleftrightarrow  
\mathbf{K}_{(c)} =\mathbf Z \mathbf A_{(c)}.
\end{align}

We note that the order of mode-$n$ products is irrelevant in a series of multiplications 
if the modes are distinct. That is,
\begin{equation}
\mathcal{X} \times_n \mathbf A \times_m \mathbf B
=
\mathcal{X} \times_m \mathbf B \times_n \mathbf A,
\qquad
\text{if } n \neq m.
\end{equation}

\subsection{Tucker Decomposition of a Third Order Tensor}

\begin{figure}[H]
\begin{minipage}{0.6\textwidth}
\centering
\includegraphics[width=\linewidth]{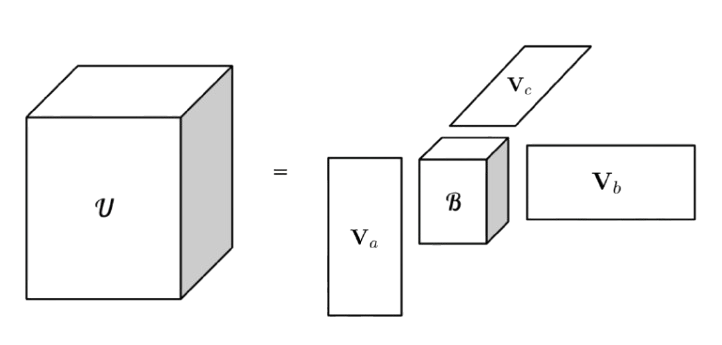}
\end{minipage}
\hfill
\begin{minipage}{0.35\textwidth}
\centering
\begin{forest}
[{$\{1,2,3\}$}
    [{$\{1\}$}]
    [{$\{2\}$}]
    [{$\{3\}$}]
]
\end{forest}
\end{minipage}
\caption{Tucker decomposition of a third order tensor ($d=3$).}
\label{fig:Tucker}
\end{figure}

Just as a matrix can be compressed using the singular value decomposition (SVD), a third order tensor can be compressed using the Tucker decomposition \cite{Tucker_1966}. Also known as the higher-order singular value decomposition (HOSVD) and higher-order PCA \cite{HOSVD}, the Tucker decomposition decomposes an $n$th-order tensor into mode-$n$ factors (view as 1D bases) and a smaller $n$th-order core tensor (view as coefficients). The Tucker decomposition for the $d=3$ case is shown in Figure \ref{fig:Tucker}. Taking the mode-$n$ product of the core tensor and the factors/bases in their respective modes reconstructs the original tensor:
\begin{equation}\label{eq:tucker_decomp}
    \mathcal{U} = \mathcal{B} \times_a \mathbf{V}_a \times_b \mathbf{V}_b \times_c  \mathbf{V}_c.
\end{equation}

If the original tensor $\mathcal{U}$ is size $N_a\times N_b \times N_c$ then the core tensor is of size $r_a\times r_b \times r_c$ and the basis matrices $\mathbf V_a,\mathbf  V_b,\ \text{and}\ \mathbf  V_c \ $are of sizes $N_a \times r_a, N_b \times r_b ,\text{and}\  N_c \times r_c$. This decomposition is said to have multilinear rank $[r_1,r_2,r_3,r_4]$. In practice, we usually represent $\mathcal{U}$ as a Tucker decomposition of low-multilinear rank, in which case the equality in \eqref{eq:tucker_decomp} is an approximation. Thus, the storage complexity has been reduced from $\mathcal{O}(N^3)$ to $\mathcal{O}(r^3+3Nr)$. If $r_a\ll N_a,r_b\ll N_b,r_c\ll N_c$, then we observe a significant reduction in the storage complexity of $\mathcal{U}$. Due to cubic scaling, these ideal low rank solutions lead to drastic increases in savings.



\subsection{Hierarchical Tucker Decomposition (HTD) of a Fourth Order Tensor}
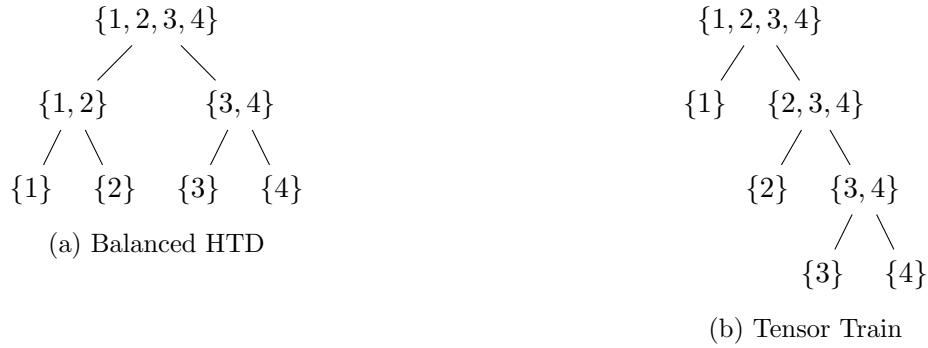
\begin{figure}[H]
\centering
\begin{subfigure}[t]{0.48\textwidth}
\centering
\vspace{0pt}
\begin{forest}
baseline=(current bounding box.north),
[{$\{1,2,3,4\}$}
    [{$\{1,2\}$}
        [{$\{1\}$}]
        [{$\{2\}$}]
    ]
    [{$\{3,4\}$}
        [{$\{3\}$}]
        [{$\{4\}$}]
    ]
]
\end{forest}
\caption{Balanced HTD}
\label{fig:balancedHTD}
\end{subfigure}
\hfill
\begin{subfigure}[t]{0.48\textwidth}
\centering
\vspace{0pt}
\begin{forest}
baseline=(current bounding box.north),
[{$\{1,2,3,4\}$}
    [{$\{1\}$}]
    [{$\{2,3,4\}$}
        [{$\{2\}$}]
        [{$\{3,4\}$}
            [{$\{3\}$}]
            [{$\{4\}$}]
        ]
    ]
]
\end{forest}
\caption{Tensor Train}
\end{subfigure}

\caption{Two examples of $d=4$ hierarchical Tucker decompositions.}
\label{fig:HTdeps}
\end{figure}

When addressing problems of high dimension $(d\geq 4)$ the $\mathcal{O}(r^d)$ storage complexity of the Tucker decomposition becomes limiting. In this paper we instead focus an extension of the Tucker decomposition known as the hierarchical Tucker Decomposition (HTD). Like the Tucker decomposition, the HTD decomposes a $d$-dimensional tensor into core and leaf tensors. But whereas the single core tensor in a Tucker decomposition is $d$-dimensional, the HTD uses a collection of third-order core tensors. In this paper, we present the $d=4$ HTD since the extension to higher dimensional HTD follows similarly, as per Remark \ref{rem:highD}. Figure \ref{fig:HTdeps} presents two possible HT decompositions of a $4$th-order tensor.

Important terminology must be considered when discussing various decompositions: The \textit{parent node}, with respect to a given node in the HT-tree diagram, refers to the node connected one level above it. Dependencies are illustrated in the diagram with solid lines. The \textit{children nodes}, with respect to a given node, are those that are connected and are one level below the given node. The \textit{leaf nodes} are on the lowest level of the tree. For example: In Figure \ref{fig:balancedHTD} the $\{1,2,3,4\}$ node is the parent node to the $\{1,2\}$ and $\{3,4\}$ nodes. And the children of the $\{1,2\}$ node are the leaf nodes $\{1\}$ and $\{2\}$. The key insight in the HTD follows from a nestedness property \cite{Kormann2017LowRankTD}. We restate it here:

\begin{thm}
\label{thm:dims}
    Let $\mathcal{U} \in \mathbb{R}^{I_1\times I_2\times ... \times I_d} \ and \ \alpha = \alpha_l\  \dot{\cup}\ \alpha_r$, where $\alpha = \{l,...,r\} \subset \{1,...,d\} $ is the parent node and the children nodes are $\alpha_l = \{l,...,m\}$ and $\alpha_r = \{{m+1},...,r\}$. Then,
\begin{equation}\label{eq:Lemma6}
    \operatorname{span}(\operatorname{mat}_{\alpha}(\mathcal{U}))\subset\operatorname{span}(\operatorname{mat}_{\alpha_r}(\mathcal{U})\otimes\operatorname{mat}_{\alpha_l}(\mathcal{U})).
\end{equation}
\end{thm}

Crucially, equation \eqref{eq:Lemma6} implies that if we have bases $\mathbf{U}_{\alpha}$, $\mathbf{U}_{\alpha_l}$, and $\mathbf{U}_{\alpha_r}$, then $\mathbf{U}_\alpha = \mathbf{B}_\alpha(\mathbf{U}_{\alpha_l}\otimes\mathbf{U}_{\alpha_r}) $, for some (flattened) core tensor $\mathbf{B}_\alpha\in \mathbb{R}^{r_{\alpha_l}r_{\alpha_r}\times r_\alpha}$.  $\mathbf{B}_{\alpha}$ can be folded back into its full tensor form $\mathcal{B}_{\alpha}\in\mathbb{R}^{r_{\alpha_l}\times r_{\alpha_r} \times r_{\alpha}}$. In summary, for a tensor $\mathcal{U}$ of arbitrary dimensions, there exists a core tensor $\mathcal{B}$ such that the $\mathcal{U}$ may be reproduced from tensors that each only depend on a subset of its dimensions. Thus, we can recursively decompose each node until arriving at leaf nodes, producing a ``tree" which eactly represents $\mathcal{U}$.

\subsubsection{The Balanced Hierarchical Tucker Decomposition in 4 Dimensions}
In this paper, we select the balanced HTD in the $d=4$ case for simplicity. Expressions for both the vectorized, and full original tensor will be shown, the latter in terms of mode-$n$ products. Equivalent decomposition equations will be presented together. 

\begin{figure}[H]
\centering
\begin{minipage}{0.15\textwidth}
\centering
\begin{forest}
   [\(\mathcal{B}_{1234}\)
    [\(\mathcal{V}_{12}\)]
    [\(\mathcal{V}_{34}\)]
  ]
\end{forest}
\end{minipage}
\hfill
\begin{minipage}{0.15\textwidth}
\centering
\begin{forest}
   [\large{\(\square\)}
    [{\includegraphics[width=.6cm]{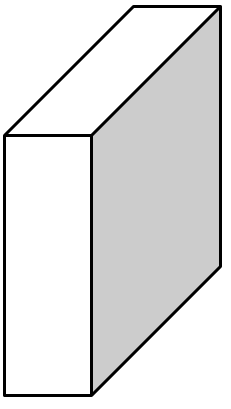}}]
    [{\includegraphics[width=.6cm]{images/U12.png}}]
  ]
\end{forest}
\end{minipage}
\hfill
\begin{minipage}{0.3\textwidth}
\centering
\begin{forest}
   [\(\mathcal{B}_{1234}\)
    [\(\mathcal{B}_{12}\)
      [\(\mathbf{V}_1\)]
      [\(\mathbf{V}_2\)]
    ]
    [\(\mathcal{B}_{34}\)
      [\(\mathbf{V}_3\)]
      [\(\mathbf{V}_4\)]
    ]
  ]
\end{forest}
\end{minipage}
\hfill
\begin{minipage}{0.3\textwidth}
\centering
\begin{forest}
   [\large{\(\square\)}
    [{\includegraphics[width=.35cm]{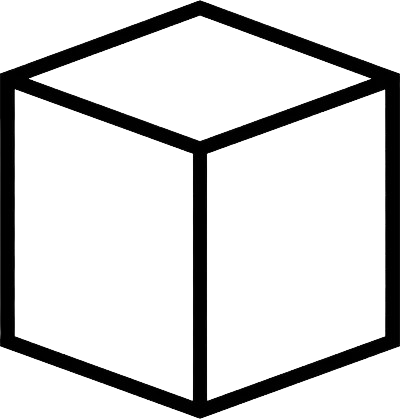}}
      [{\includegraphics[width=.5cm]{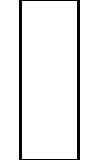}}]
      [{\includegraphics[width=.5cm]{images/rectangle.png}}]
    ]
    [{\includegraphics[width=.35cm]{images/cube.png}}
      [{\includegraphics[width=.5cm]{images/rectangle.png}}]
      [{\includegraphics[width=.5cm]{images/rectangle.png}}]
    ]
  ]
\end{forest}
\end{minipage}
\caption{Hierarchical Tucker Decomposition in 4D}
\label{fig:HTtree}
\end{figure}

Due to its prominence in our algorithm, for arbitrary tensors, $\mathcal{B}$, bold face $\mathbf{B}$ will be used to indicate a transpose of the flattening of the tensor along the parent mode ($\text{mat}_c({\mathcal{B}})^T$). We note that although the root transfer tensor $\mathbf{B}_{1234}$ is a $r_{12}r_{34}\times 1$ column vector, for notational convenience we instead define $\operatorname{vec}(\mathbf{B}_{1234})$ as the root transfer tensor/matrix, with $\mathbf{B}_{1234}\in\mathbb{R}^{r_{12}\times r_{34}}$ as its reshaped matrix. By equation \eqref{eq:Lemma6}, we can begin the balanced binary decomposition of a four-way tensor in the following fashion. 
\begin{equation}\label{eq:htd4_-1}
    \text{vec}(\mathcal{U}) = (\mathbf{V}_{34} \otimes \mathbf{V}_{12} )\operatorname{vec}(\mathbf{B} _{1234})
\end{equation}

We note that this is equivalent to 
\begin{equation}\label{eq:htd4_-1}
    \mathcal{U} = \mathcal{B}_{1234}\times_{a}\mathbf{ V}_{12}\times_{b}\mathbf{V}_{34}
\end{equation}

To further decompose these tensors, we may perform subsequent Tucker decompositions at each parent node. For our 4d tensor, at the $\{1,2\}$ and $\{3,4\}$ nodes,
\begin{subequations}
\begin{align}
\mathcal{V}_{34} = \mathcal{B}_{34} \times_a\mathbf{V}_3\times_b\mathbf V_4,& \qquad
\mathcal{V}_{12} = \mathcal{B}_{12} \times_a\mathbf{V}_1\times_b\mathbf V_2,
\label{eq:htd4_0_a}\\
\mathbf{V}_{34} = (\mathbf{V}_4 \otimes \mathbf{V}_3)\, \mathbf{B}_{34},& \qquad
\mathbf{V}_{12} = (\mathbf{V}_2 \otimes \mathbf{V}_1)\, \mathbf{B}_{34}.
\label{eq:htd4_0_b}
\end{align}
\end{subequations}

We can see that our mode-$n$ product formulation of $\eqref{eq:htd4_0_b}$ is equivalent to decomposing $\mathcal{V}_{12}$ and $\mathcal{V}_{34}$ via the Tucker decomposition. Rearranging \eqref{eq:htd4_-1} using Kronecker Product \textit{Property 1} yields the final balanced binary HT decomposition of a fourth-order tensor.
\begin{align}\label{eq:htd4}
    \text{vec}\big(\mathcal{U}\big) = (\mathbf{V}_4 \otimes \mathbf{V}_3 \otimes \mathbf{V}_2 \otimes \mathbf{V}_1)\,
    (\mathbf{B}_{34} \otimes \mathbf{B}_{12})\, \operatorname{vec}(\mathbf{B}_{1234})
\end{align}

The size of each basis $\mathbf{V}_i$ is the mode length times the rank in that dimension, which is generally much smaller by comparison. ($N_i\times r_i, \ r_i\ll N_i$ ), and the size of the core tensors scales with the product of the ranks of its leaf and parent tensors. We see that, for example $\mathcal{B}_{12}$ has size $r_1\times r_2\times r_{12}$, and $\mathbf{B}_{12}$ has size $r_1r_2\times r_{12} $ as it is flattened in the parent node and transposed. We compute that this decomposition has a storage complexity given by 
\begin{equation}\label{eq:HT4complx}
    \mathcal{O}(4Nr+2r^3+r^2)
\end{equation}
For a $d$-dimensional tensor, the HT decomposition exhibits a storage complexity given by
\begin{equation} \label{HTdcomplx}
    \mathcal{O}(dNr+(d-2)r^3+r^2)
\end{equation}

We recall that the complexity for the Tucker Decomposition is $\mathcal{O}(dNr+r^d)$, in comparison to the complexity for the HTD at $\mathcal{O}(dNr+dr^3)$. We see that for $d=3$, the HTD and Tucker decomposition have comparable complexity, but for $d\geq 4$ the HTD becomes advantageous. 

\begin{rem}
We can view each core tensor as being the core of a third-order Tucker tensor local to that node, connecting the parent basis and the children bases. For example: $\mathcal{U} = \mathcal{B}_{12}\times_a\mathbf {V}_1\times_b \mathbf{V}_2 \times_c \mathbf{V}_{12}$.
\end{rem}


\section{Algorithm for 4-dimensional Diffusion Equations}

For simplicity, we consider the $d=4$ case with the $d>4$ cases described in a later remark. We now describe the implicit scheme for solving the diffusion equation
\begin{equation}\label{eq:4dheat}
    u_t =d_1(t)u_{x_1x_1}+ d_2(t)u_{x_2x_2}+ d_3(t)u_{x_3x_3}+ d_4(t)u_{x_4x_4}.
\end{equation}

Naively storing the entire solution in a fourth-order tensor is too expensive, so we utilize the HTD. 
After storing the HT solution, we must update the core tensors and leaf bases separately. We work from the bases upward, in reference to the structure presented in figure \ref{fig:HTtree}, to the core tensors, passing through each subtree. 
Discretizing equation \eqref{eq:4dheat} \textit{in space} over uniform computational grids presented in \eqref{eq:compgrid},
we assume that the vectorized form of the 4D numerical solution can be represented as a time-dependent fourth-order tensor in Hierarchical Tucker format, beginning with
\begin{equation}\label{eq:htd42}
    \operatorname{vec}(\mathcal{U}(t))
    = \big( \mathbf{V}_{34}(t) \otimes \mathbf{V}_{12}(t) \big)\, \operatorname{vec}(\mathbf{B}_{1234}(t)).
\end{equation}
We consider an expression for the time-dependent solution $\mathcal{U}(t)$ as refoldings of \eqref{eq:htd42}, which by Property \ref{prop2} is 
\begin{align}
\label{eq:Btil_0_12}
    \mathbf{V}_{12}(t)\, \mathbf{B}_{1234}(t) \mathbf{V}_{34}^T(t).
\end{align}

Without loss of generality, we will now focus on the left subtree of the tree tensor network presented in Figure \ref{fig:HTtree}. By equation \eqref{eq:htd4_0_b}, we can further expand $\mathbf{V}_{12}$ to get
\begin{align}\label{eq:Btil_step0_12}
    (\mathbf{V}_{2}(t)\otimes \mathbf V_{1}(t)) \mathbf{B}_{12}(t) \, \mathbf{B}_{1234}(t) \mathbf{V}^T_{34}(t).
\end{align}

Defining the projected core tensor, $\tilde{\mathcal{B}}_{12}\in \mathbb{R}^{N_1\times N_2 \times r_{34}}$  where $\tilde{\mathcal{B}}_{12}\coloneq\mathcal{B}_{12}\times_c\mathbf{B}_{1234}$ such that $\tilde{\mathbf B}_{12} = \mathbf{B}_{12} \, \mathbf{B}_{1234} $, we can refold solution  \eqref{eq:Btil_step0_12} using mode-$n$ products to get
\begin{align} \label{U_fold_12}
    \tilde{\mathcal{B}}_{12}(t)\times_a \mathbf{V}_1(t)\times_b\mathbf{V}_2(t)\times_c\mathbf{V}_{34}(t).
\end{align}

We note the resemblance of our formulations in equation \eqref{U_fold_12} to a third-order Tucker decomposition and may now build off of the work of other projection-based integrators that utilize the Tucker Decomposition. Discretizing equation \eqref{eq:4dheat} in space under the decomposition \eqref{U_fold_12} (local to this subtree) yields the tensor differential equation
\begin{align}
     \begin{split}
        \frac{d}{dt} \Big[\tilde{\mathcal{B}}_{12}(t)\times_a \mathbf{V}_1(t)\times_b\mathbf{V}_2(t)\times_c&\mathbf{V}_{34}(t)\Big] =\tilde{\mathcal{B}}_{12}(t)\times_a \mathbf F_1\mathbf{V}_1(t)\times_b\mathbf{V}_2(t)\times_c\mathbf{V}_{34}(t)\\
        &+\tilde{\mathcal{B}}_{12}(t)\times_a \mathbf{V}_1(t)\times_b\mathbf F_2\mathbf{V}_2(t)\times_c\mathbf{V}_{12}(t)\\
        &+\tilde{\mathcal{B}}_{12}(t)\times_a \mathbf{V}_1(t)\times_b\mathbf{V}_2(t)\times_c(\mathbf{I}_{N_4}\otimes \mathbf{F}_3\mathbf{V}_3(t))\mathbf{V}_{34}(t)\\
        &+\tilde{\mathcal{B}}_{12}(t)\times_a \mathbf{V}_1(t)\times_b\mathbf{V}_2(t)\times_c(\mathbf{F}_{4}\mathbf V_4(t)\otimes \mathbf{I}_{N_3})\mathbf{V}_{34}(t),
        \label{eq:ddt}
    \end{split}
\end{align}
where for $i\in \{1,2,3,4\}$, $\mathbf F_i$ represents the discretization of the one dimensional Laplacian $d_i\partial_{x_i}^2$, e.g., from a finite difference or spectral method. 


\subsection{The first-order scheme using implicit Euler}\label{subsec:firstorder}
In this section, we will specifically examine the left side of the tree tensor network presented in figure \ref{fig:HTtree}. 
First, we will update the leaf bases. Without loss of generality, we will explicitly show only $\mathbf V_1$ and $\mathbf{V}_2$, as $\mathbf{V}_3$ and $\mathbf{V}_4$ follow similarly. This process is called the \textit{K-steps}. We will then proceed up the tree to the core tensor $\mathcal{B}_{12}$ and $\mathcal{B}_{34}$ and finally $\mathcal{B}_{1234}$ during the \textit{B-steps}. 
\subsubsection*{K-steps}
Integrating equation \eqref{eq:ddt} using the implicit Euler method yields
\begin{align}
\begin{split}
        \tilde{\mathcal{B}}_{12} ^{n+1}\times_a\mathbf{V}_1 ^{n+1}\times_b\mathbf{V}_2 ^{n+1}\times_c&\mathbf{V}_{34} ^{n+1} 
        =\tilde{\mathcal{B}}_{12} ^{n}\times_a \mathbf{V}_1^{n}\times_b\mathbf{V}_2^{n}\times_c\mathbf{V}_{34}^{n}\\
        &+\Delta t\Big\{\tilde{\mathcal{B}}_{12} ^{n+1}\times_a \mathbf F_1\mathbf{V}_1^{n+1}\times_b\mathbf{V}_2^{n+1}\times_c\mathbf{V}_{34}^{n+1}\\
        &+\tilde{\mathcal{B}}_{12}^{n+1}\times_a \mathbf{V}_1^{n+1}\times_b\mathbf F_2\mathbf{V}_2^{n+1}\times_c\mathbf{V}_{34}^{n+1}\\
        &+\tilde{\mathcal{B}}_{12}^{n+1}\times_a \mathbf{V}_1^{n+1}\times_b\mathbf{V}_2^{n+1}\times_c(\mathbf{V}^{n+1}_{4}\otimes \mathbf{F}_3\mathbf{V}_3^{n+1})\mathbf{B}_{34}^{n+1}\\
        &+\tilde{\mathcal{B}}_{12}^{n+1}\times_a \mathbf{V}_1^{n+1}\times_b\mathbf{V}_2^{n+1}\times_c(\mathbf{F}_{4}\mathbf{V}_4^{n+1}\otimes \mathbf{V}^{n+1}_{3})\mathbf{B}_{34}^{n+1} \Big\}.
        \label{eq:impup}
\end{split}
\end{align}

In order to update each leaf basis, we project equation \eqref{eq:impup} onto subspaces that span all but a single dimension. In order to project equation \eqref{eq:impup}, we consider $\mathbf{V}_i^{\star,{n+1}} $ such that $\mathbf{V}_i(t^{n+1}) \approx\mathbf{V}_i^{\star,{n+1}}$ since we do not yet know the future bases. For the first order accurate implicit Euler method, we let $\mathbf{V}_i^{\star,{n+1}} = \mathbf{V}_i^n$. The projected solutions are
\begin{equation}\label{eq:K1_K2}
\begin{aligned}
\mathbf{K}_{1}^{n+1}
& \coloneq  
\mathbf{V}_{1}^{n+1}\,
\operatorname{mat}_a\!\Big( \tilde{\mathcal{B}}_{12}^{n+1}\Big)
\Big(
(\mathbf{V}_{34}^{n+1})^T\mathbf{V}_{34}^{\star,n+1}
\;\otimes\;
(\mathbf{V}_{2}^{n+1})^T\mathbf{V}_{2}^{\star,n+1}
\Big),
\\
\mathbf{K}_{2}^{n+1}
&\coloneq
\mathbf{V}_{2}^{n+1}\,
\operatorname{mat}_b\!\Big( \tilde{\mathcal{B}}_{12}^{n+1}\Big)
\Big(
(\mathbf{V}_{34}^{n+1})^T\mathbf{V}_{34}^{\star,n+1}
\;\otimes\;
(\mathbf{V}_{1}^{n+1})^T\mathbf{V}_{1}^{\star,n+1}
\Big),
\end{aligned}
\end{equation}
where $\mathbf{V}_{34} = \big(\mathbf{V}_4\otimes\mathbf{V}_3 \big)\mathbf{B}_{34}$. To update the leaf bases in modes 1 and 2, we flatten equation \eqref{eq:impup} in the modes $a$ or $b$, respectively. For example since $\mathbf{K}_{1} $ no longer has dependence on modes $2,3, \text{and} \ 4$, we can use it to extract an orthonormal basis in mode $1$. Without loss of generality, we will present only the update equation for $\mathbf V_1$. Performing the mode $a$ matricization on equation \eqref{eq:impup} we obtain
\begin{align}\label{eq:impup_mat}
\begin{split}
    \mathbf{V}_1^{n+1}\tilde{\mathbf{B}}_{12,(a)}^{n+1}&\big(\mathbf{V}_{34}^{n+1} \otimes \mathbf{V}_2^{n+1}\big)^{T} = \mathbf{V}_1^{n}\tilde{\mathbf{B}}_{12,(a)}^{n}\big(\mathbf{V}_{34}^{n} \otimes \mathbf{V}_2^{n}\big)^{T}\\
    &+ \Delta t\Big\{ \mathbf{F}_1 \mathbf{V}_1^{n+1}\tilde{\mathbf{B}}_{12,(a)}^{n+1}\big(\mathbf{V}_{34}^{n+1} \otimes \mathbf{V}_2^{n+1}\big)^{T} + \mathbf{V}_1^{n+1}\tilde{\mathbf{B}}_{12,(a)}^{n+1}\big(\mathbf{V}_{34}^{n+1} \otimes \mathbf{F}_2 \mathbf{V}_2^{n+1}\big)^{T}\\
    &+ \mathbf{V}_1^{n+1}\tilde{\mathbf{B}}_{12,(a)}^{n+1}\Big((\mathbf{V}_{4}^{n+1} \otimes \mathbf{F}_3 \mathbf{V}_3^{n+1})\mathbf{B}_{34}^{n+1}\otimes \mathbf{V}_2^{n+1}\Big)^{T}\\
    &+ \mathbf{V}_1^{n+1}\tilde{\mathbf{B}}_{12,(a)}^{n+1}\Big((\mathbf{F}_4 \mathbf{V}_4^{n+1}\otimes\mathbf{V}_{3}^{n+1})\mathbf{B}_{34}^{n+1}\otimes \mathbf{V}_2^{n+1}\Big)^{T}\Big\}.
\end{split}
\end{align}

We then project equation \eqref{eq:impup_mat} onto $\mathbf{V}_2^{\star,n+1} \coloneq \mathbf V_2^n$ and $\mathbf{V}_{34}^{\star,n+1} \coloneq\big(\mathbf{V}_4^n\otimes\mathbf{V}_3^n \big)\mathbf{B}_{34}^n$. Multiplying equation \eqref{eq:impup_mat} on the right by $\mathbf{V}_{34}^{\star,n+1}\otimes\mathbf{V}_2^{\star,n+1}$, and making use of Properties \ref{prop1}-\ref{prop2} and the definition of $\mathbf{K}_1$, we obtain the Sylvester equation
\begin{equation}\label{eq:sylv}
    \mathbf{A}\mathbf{K}_1^{n+1} + \mathbf{K}_1^{n+1}\mathbf{B} = \mathbf{C},
\end{equation}
where $\mathbf{A} = \mathbf{I}_{N_1}-\Delta t\mathbf{F}_1$, $\mathbf{C}=\mathbf{K}_1^n$, and
\begin{align}
\begin{split}
    \mathbf{B} =& - \Delta t\Big( (\mathbf{B}_{34}^{\star,n+1})^T((\mathbf{F}_4\mathbf{V}_4^{\star,n+1})^T\mathbf{V}_4^{\star,n+1}\otimes\mathbf{I}_{r_3})\mathbf{B}_{34}^{\star,n+1}\otimes\mathbf{I}_{r_2}\\
    &+ (\mathbf{B}_{34}^{\star,n+1})^T(\mathbf{I}_{r_4}\otimes(\mathbf{F}_3\mathbf{V}_3^{\star,n+1})^T\mathbf{V}_3^{\star,n+1})\mathbf{B}_{34}^{\star,n+1}\otimes\mathbf{I}_{r_2}\\
    &+(\mathbf{I}_{r_{34}}\otimes (\mathbf{F}_2\mathbf{V}_2^{\star,n+1})^T\mathbf{V}_2^{\star,n+1}
     \Big).
\end{split}
\end{align}

We solve equation \eqref{eq:sylv} using the built-in \texttt{sylvester} function from MATLAB, although any appropriate Sylvester solver can be used. We can now extract the orthonormal basis $\mathbf{V}^{n+1}_1$ with a reduced QR factorization, $\mathbf{K}^{n+1}_1 = \mathbf{Q}\mathbf{R} =: \mathbf{V}_1^{\ddagger,n+1}\mathbf{R}$. Similarly one can obtain and solve the Sylvester equations for $\mathbf{K}_i$, $ i = 2,3,4$. Here, the double dagger $\ddagger$ denotes the one-dimensional orthonormal bases obtained from the K-steps.

We adopt the \textit{reduced augmentation procedure} from \cite{2d_rail} and augment each $\mathbf{V}^{\ddagger,n+1}_i$ with the previous basis $\mathbf{V}^n_i$ for $i = 1,2,3,4$. Augmenting the updated basis with the current one is advantageous because it allows the spanning subspace to retain information over the full interval $[t^n, t^{n+1}]$. This becomes particularly important when the solution changes rapidly, for example, during short-time diffusive dynamics. We compute the reduced QR factorizations of the augmented bases as follows: 
\begin{equation}
\left[\mathbf{V}_i^{\ddagger,n+1}, \mathbf{V}_i^n \right]
= \mathbf{Q}_i \mathbf{R}_i, \ i = 1,2,3,4.
\end{equation}

To avoid unnecessarily doubling the rank, we compress/reduce our augmented bases. We first compute the singular value decomposition (SVD) of $\mathbf{R}$, that is, $\mathbf{R}_i=\mathbf{U}_i\mathbf{\Sigma}_i\mathbf{V}_i$. In each mode, let $\hat{\mathbf{U}}_i$ be the left singular vectors of $\mathbf{R}_i$ that correspond to singular values greater than a fixed tolerance $10^{-12}$. This tolerance is large enough to remove redundant information but small enough to have a negligible effect on the scheme's consistency. 
We let $\hat{r}_i$ denote the rank of $\hat{\mathbf{U}}_i$. In each mode, we define the reduced augmented basis $\hat{\mathbf{V}}_i^{n+1}\coloneq \mathbf{Q}_i\hat{\mathbf{U}}_i$. Now that we have computed $\hat{\mathbf{V}}_i^{n+1}$ for $i=1,2,3,4$ we may proceed to updating the core tensors. 

\subsubsection*{B-steps}

Now that we have updated the one-dimensional bases, we can proceed to update the core tensors. Without loss of generality, we consider $\mathcal{B}_{12}$ since $\mathcal{B}_{34}$ follows similarly. The core tensor at the $\{1,2\}$ node can now be updated via a Galerkin projection. We vectorize equation \eqref{eq:impup} and obtain 
\begin{align}
\label{eq:impup_vec}
    \begin{split} 
        (\mathbf{V}_{34}^{n+1} \otimes \mathbf{V}_{2}^{n+1}&\otimes \mathbf{V}_1^{n+1})\operatorname{vec}(\tilde{\mathcal{B}}_{12}^{n+1}) = (\mathbf{V}_{34}^{n} \otimes \mathbf{V}_{2}^{n}\otimes \mathbf{V}_{1}^{n})\operatorname{vec}(\tilde{\mathcal{B}}_{12}^{n})\\
        &+ \Delta t \Big\{ (\mathbf{V}_{34}^{n+1} \otimes \mathbf{V}_{2}^{n+1}\otimes \mathbf{F}_1\mathbf{V}_1^{n+1}) + (\mathbf{V}_{34}^{n+1}\otimes\mathbf{F}_2\mathbf{V}_{2}^{n+1}\otimes \mathbf{V}_1^{n+1})\\
        & +\big((\mathbf{V}_{4}^{n+1} \otimes \mathbf{F}_3\mathbf{V}_3^{n+1})\mathbf{B}_{34}^{n+1} \otimes \mathbf{V}_{2}^{n+1}\otimes \mathbf{V}_1^{n+1}\big)\\
        & +\big((\mathbf{F}_4\mathbf{V}_{4}^{n+1} \otimes \mathbf{V}_3^{n+1})\mathbf{B}_{34}^{n+1} \otimes \mathbf{V}_{2}^{n+1}\otimes \mathbf{V}_1^{n+1}\big)\Big\}\operatorname{vec}(\tilde{\mathcal{B}}_{12}^{n+1}).
    \end{split}
\end{align}

We now perform a Galerkin projection on equation \eqref{eq:impup_vec} using the updated leaf bases $\hat{\mathbf{V}}_1^{n+1}$ and $\hat{\mathbf{V}}_2^{n+1}$, and $\mathbf{V}_{34}^{\star,n+1}$. We note that we use $\mathbf{V}^{\star,n+1}_{34}$ because we have not yet updated the core tensor at the $\{3,4\}$ node. We define \begin{equation}\label{eq:vecBtilhat}
\operatorname{vec}(\hat{\tilde{\mathcal B}}^{n+1}_{12} )\ \coloneq \left((\mathbf{V}_{34}^{\star,n+1})^{T} \mathbf{V}_{34}^{{n+1}}
    \otimes (\hat{\mathbf{V}}_{2}^{n+1})^T \mathbf{V}_{2}^{{n+1}}
    \otimes (\hat{\mathbf{V}}_{1}^{n+1})^T\mathbf{V}_{1}^{{n+1}}
  \right)
\operatorname{vec}(\tilde{\mathcal{B}}_{12}^{{n+1}}).
\end{equation}

Multiplying equation \eqref{eq:impup_vec} on the left by $(\mathbf{V}_{34}^{\star,n+1}\otimes\hat{\mathbf{V}}_2^{n+1}\otimes\hat{\mathbf{V}}_1^{n+1})^T$, after some tedious but straightforward algebra we obtain the third order tensor linear equation

\begin{align}\label{eq:impup_vec_proj_lin}
\begin{split}
\Bigg\{&
(\mathbf{B}_{34}^{\star,n+1})^T
\Big(
  -\Delta t (\hat{\mathbf{V}}_4^{n+1})^T
  (\mathbf{F}_4 \hat{\mathbf{V}}_4^{n+1})
  \otimes \mathbf{I}_{\hat{r}_3^{n+1}}
\Big)
\mathbf{B}_{34}^{\star,n+1} \otimes \mathbf{I}_{\hat{r}_2^{n+1}} \otimes \mathbf{I}_{\hat{r}_1^{n+1}}\\
&\qquad+ (\mathbf{B}_{34}^{\star,n+1})^T
\Big(
  \mathbf{I}_{\hat{r}_4^{n+1}}
  \otimes
  -\Delta t (\hat{\mathbf{V}}_4^{n+1})^T
  (\mathbf{F}_4 \hat{\mathbf{V}}_4^{n+1})
\Big)
\mathbf{B}_{34}^{\star,n+1} \otimes \mathbf{I}_{\hat{r}_2^{n+1}} \otimes \mathbf{I}_{\hat{r}_1^{n+1}}\\
&\qquad\qquad+ \mathbf{I}_{\hat{r}_{34}^{n+1}}
\otimes
\Big(
  -\Delta t (\hat{\mathbf{V}}_2^{n+1})^T
  (\mathbf{F}_2 \hat{\mathbf{V}}_2^{n+1})
\Big)
\otimes \mathbf{I}_{\hat{r}_1^{n+1}}\\
&\qquad\qquad\quad+ \mathbf{I}_{\hat{r}_{34}^{n+1}}
\otimes \mathbf{I}_{\hat{r}_2^{n+1}}
\otimes
\Big(
  -\Delta t (\hat{\mathbf{V}}_1^{n+1})^T
  (\mathbf{F}_1 \hat{\mathbf{V}}_1^{n+1})
\Big)
\Bigg\}
\operatorname{vec}(\hat{\tilde{\mathcal{B}}}_{12}^{n+1})
= \operatorname{vec}(\hat{\tilde{\mathcal{B}}}_{12}^{n}).
\end{split}
\end{align}


Naively solving equation \eqref{eq:impup_vec_proj_lin} by taking the inverse can be very expensive. Following the work in \cite{2d_rail}, we use the direct solver proposed by Simoncini in \cite{Simoncini2020} extended to a general Tucker tensor $\mathcal{B}$. Simoncini presents the algorithm for the rank-1 case, but the extension to a general $\mathcal{B}$ is straightforward and was shown in \cite{2d_rail}. For brevity, we refer the reader to \cite{Simoncini2020} for more details of the direct solver. We proceed having computed $\operatorname{vec}(\hat{\tilde{\mathcal{B}}}_{12}^{n+1})$.

To extract $\mathcal{B}_{12}^{n+1}$ we reshape $\operatorname{vec}(\hat{\tilde{\mathcal{B}}}_{12}^{n+1}) $ into $\hat{\tilde{\mathbf{B}}}_{12}^{n+1}\in \mathbb{R}^{\hat{r}^{n+1}_1\hat{r}^{n+1}_2\times\hat{r}^{n+1}_{34}}$. Recalling that $\tilde{\mathbf{B}}_{12} = {\mathbf{B}}_{12}{\mathbf{B}}_{1234} $, we see that $\mathbf{B}_{12}$ and $\tilde{\mathbf{B}}_{12}$ share the same column space. As such, we can now extract the $\hat{\tilde{\mathbf{B}}}_{12}^{n+1}$ by taking a reduced QR factorization $\hat{\tilde{\mathbf{B}}}_{12}^{n+1}=\mathbf{Q}_{12}\mathbf{R}_{12}$. We let our flattened updated core tensor be $\hat{\mathbf{B}}_{12}^{n+1} \coloneqq \mathbf{Q}_{12}$.
We have now computed the updated leaf bases $\hat{\mathbf{V}}_{1}^{n+1}$, $\hat{\mathbf{V}}_{2}^{n+1}$ and core tensor $\hat{\mathcal{B}}_{12}^{n+1}$. We can follow the same procedures analogously to determine $\hat{\mathbf{V}}_{3}^{n+1}$, $\hat{\mathbf{V}}_{4}^{n+1}$ and $\hat{\mathcal{B}}_{34}^{n+1}$ for the right subtree. 

Finally, to update $\mathcal{B}_{1234}$ we look back to the flattening of $\mathcal{U}(t)$ stated in equation \eqref{eq:Btil_0_12}, and project out all spatial dependence by projecting onto $\hat{\mathbf{V}}_{12}^{n+1} = (\hat{\mathbf{V}}_2^{n+1}\otimes\hat{\mathbf{V}}_1^{n+1})\hat{\mathbf{B}}^{n+1}_{12}\ \text{and}\ \hat{\mathbf{V}}_{34}^{n+1} = (\hat{\mathbf{V}}_4^{n+1}\otimes\hat{\mathbf{V}}_3^{n+1})\hat{\mathbf{B}}^{n+1}_{34}$. Discretizing equation \eqref{eq:4dheat} local to the root node with implicit Euler, we obtain the following equation: 

\begin{align}
\begin{split}
\label{eq:b1234}
\mathbf{V}_{12}^{n+1}\;
\mathbf{B}_{1234}^{n+1}\;
(\mathbf{V}_{34}^{n+1})^T
=&
\mathbf{V}_{12}^{n}\;
\mathbf{B}_{1234}^{n}\;
(\mathbf{V}_{34}^{n})^T\\
&+\Delta t \Big\{
(\mathbf{V}_2^{n+1}\otimes \mathbf{F}_1\mathbf{V}_{1}^{n+1})
{{\mathbf B}}_{12}^{n+1}\;
\mathbf{B}_{1234}^{n+1}\;
(\mathbf{V}_{34}^{n+1})^T\\
&+
(\mathbf{F}_2\mathbf{V}_2^{n+1}\otimes \mathbf{V}_{1}^{n+1})
{{\mathbf B}}_{12}^{n+1}\;
\mathbf{B}_{1234}^{n+1}\;
(\mathbf{V}_{34}^{n+1})^T
\\
&+\mathbf{V}_{12}^{n+1}\;
\mathbf{B}_{1234}^{n+1}\;
(\mathbf{B}_{34}^{n+1})^T
(\mathbf{V}_4^{n+1}\otimes \mathbf{F}_3\mathbf{V}_3^{n+1})^T\\
&+\mathbf{V}_{12}^{n+1}\;
\mathbf{B}_{1234}^{n+1}\;
(\mathbf{B}_{34}^{n+1})^T
(\mathbf{F}_4\mathbf{V}_4^{n+1}\otimes \mathbf{V}_3^{n+1})^T
\Big\}.
\end{split}
\end{align}
We now define 
\begin{equation}\label{eq:Bhat1234def}
    \hat{\mathbf{B}}^{n+1}_{1234}\coloneq (\hat{\mathbf{V}}_{12}^{n+1})^T\mathbf{V}_{12}^{n+1} \mathbf{B}_{1234}^{{n+1}}(\mathbf{V}_{34}^{n+1})^T\hat{\mathbf{V}}_{34}^{n+1}
\end{equation}
and project equation \eqref{eq:b1234} onto $\hat{\mathbf{V}}_{12}^{n+1}$ and $\hat{\mathbf{V}}_{34}^{n+1}$ to obtain the Sylvester equation 
\begin{align}\label{eq:b1234_sylv}
\begin{split}
\Big(
\mathbf{I}_{\hat{r}_{12}}
- \Delta t &\, (\hat{\mathbf B}_{12}^{n+1})^T
\big(
(\hat{\mathbf V}_2^{n+1})^T \mathbf{F}_2 \mathbf{V}_2^{n+1}
\otimes (\hat{\mathbf V}_1^{n+1})^T \mathbf{V}_1^{n+1}\big)
\\
&+ (\hat{\mathbf V}_2^{n+1})^T \mathbf{V}_2^{n+1}
\otimes (\hat{\mathbf V}_1^{n+1})^T \mathbf{F}_1 \mathbf{V}_1^{n+1}
\big)
\mathbf{B}_{12}^{n+1}
\Big)
\hat{\mathbf{B}}_{1234}^{n+1}
\\
&\qquad+ \hat{\mathbf{B}}_{1234}^{n+1}
\Big(
- \Delta t \, (\hat{\mathbf B}_{34}^{n+1})^T
\big(
(\mathbf{V}_4^{n+1})^T \hat{\mathbf V}_4^{n+1}
\otimes (\mathbf{F}_3 \mathbf{V}_3^{n+1})^T \hat{\mathbf V}_3^{n+1}
\\
&\qquad\qquad+ (\mathbf{F}_4 \mathbf{V}_4^{n+1})^T \hat{\mathbf V}_4^{n+1}
\otimes (\mathbf{V}_3^{n+1})^T \hat{\mathbf V}_3^{n+1}
\big)
\hat{\mathbf B}_{34}^{n+1}
\Big)
=
\hat{\mathbf{B}}_{1234}^{n}.
\end{split}
\end{align}

With all of our components updated (denoted with hats), we utilize the standard hierarchical Tucker Tensor truncation algorithm included in the HT toolbox \cite{10.1145/2538688} to truncate our updated HT tensor solution according to a tolerance $\varepsilon = 10^{-6}$. We now have the final updated and truncated solution $\text{vec}\big(\mathcal{U}^{n+1}\big) = (\mathbf{V}^{n+1}_4 \otimes \mathbf{V}^{n+1}_3 \otimes \mathbf{V}^{n+1}_2 \otimes \mathbf{V}^{n+1}_1)\,
    (\mathbf{B}^{n+1}_{34} \otimes \mathbf{B}^{n+1}_{12})\, \operatorname{vec}(\mathbf{B}^{n+1}_{1234})$.

\subsection{High Order Schemes using Diagonally Implicit RK-methods}

We extend the first-order implicit scheme to high-order accuracy using diagonally implicit Runge-Kutta (DIRK) methods \cite{ascher1997implicit}. A general $s$ stage DIRK method is typically expressed in a Butcher tableau given by Table~\ref{tab:dirk}, where $c_k = \sum_{j=1}^{k} a_{kj}$ for $k = 1,2,\ldots,s$, and
$\sum_{k=1}^{s} b_k = 1$ for consistency.
Each row in the Butcher tableau represents an intermediate stage for the
solution at time $t^{(k)} = t^n + c_k \Delta t$, for  $k = 1,2,\ldots,s$.

\begin{table}[h]
\centering
\[
\begin{array}{c|cccc}
c_1 & a_{11} & 0      & \cdots & 0 \\
c_2 & a_{21} & a_{22} & \cdots & 0 \\
\vdots & \vdots & \vdots & \ddots & \vdots \\
c_s & a_{s1} & a_{s2} & \cdots & a_{ss} \\
\hline
    & b_1 & b_2 & \cdots & b_s
\end{array}
\]
\caption{Butcher tableau for an $s$ stage DIRK scheme}
\label{tab:dirk}
\end{table}

We follow the same general procedure as the first-order method (K-steps, B-steps, truncation) at each stage of the DIRK method. However, at each stage, we enrich the subspaces that we project onto with the information obtained at the previous stages. It suffices to present the high-order extension of the equations presented in Section \ref{subsec:firstorder}. In particular, the intermediate equation updating the left side of the tree at the $k$th stage is
\begin{equation}
     \begin{aligned}
        \tilde{\mathcal{B}}_{12} ^{(k)}\times_a \mathbf{V}_1 ^{(k)}\times_b\mathbf{V}_2 ^{(k)}\times_c\mathbf{V}_{34} ^{(k)} 
        =&\tilde{\mathcal{B}}_{12} ^{n}\times_a \mathbf{V}_1^{n}\times_b\mathbf{V}_2^{n}\times_c\mathbf{V}_{34}^{n}\\ +\Delta t \sum_{l=1}^{k}a_{k\ell}\Big\{&\tilde{\mathcal{B}}_{12} ^{(\ell)}\times_a \mathbf F_1\mathbf{V}_1^{(\ell)}\times_b\mathbf{V}_2^{(\ell)}\times_c\mathbf{V}_{34}^{(\ell)}
        \\+&\tilde{\mathcal{B}}_{12}^{(\ell)}\times_a \mathbf{V}_1^{(\ell)}\times_b\mathbf F_2\mathbf{V}_2^{(\ell)}\times_c\mathbf{V}_{34}^{(\ell)}\\
        +&\tilde{\mathcal{B}}_{12}^{(\ell)}\times_a \mathbf{V}_1^{(\ell)}\times_b\mathbf{V}_2^{(\ell)}\times_c(\mathbf{V}_{4}\otimes \mathbf{F}_3\mathbf{V}_3^{(\ell)})\mathbf{B}_{34}^{(\ell)}
        \\
        +&\tilde{\mathcal{B}}_{12}^{(\ell)}\times_a \mathbf{V}_1^{(\ell)}\times_b\mathbf{V}_2^{(\ell)}\times_c(\mathbf{F}_{4}\mathbf{V}_3^{(\ell)}\otimes \mathbf{V}_{3})\mathbf{B}_{34}^{(\ell)} \Big\}.
    \end{aligned} \label{eq:impuprk}
\end{equation}

We consider stiffly accurate DIRK methods where $c_s =1$ and $a_{sk} = b_k$ such that $t^{(s)} = t^{n+1}$. The Butcher tableaus for the second and third-order DIRK methods are included in the Appendix. Unlike first first-order scheme where we projected onto the current basis $t^n$, the approximate basis is enriched with information from the previous stages. To enrich the approximate bases we perform the reduction augmentation concept on the solutions from each time step to obtain a solution in HT format. Since the core tensors are included in the projections, augmenting the 1D bases alone is insufficient. We sum HT solutions from previous stages, then compress/truncate (tolerance = $1e-6$). 
\begin{equation}\label{eq:redaug_HT}
    HT(\mathcal{U}^{k,\star}) =  \texttt{redaug}\left(HT(\mathcal{U}^{\dagger,(k)}), HT(\mathcal{U}^{(k-1)}),..., HT(\mathcal{U}^{n})\right)
\end{equation}
where $HT(\mathcal{U}^{\dagger,(k)})$ denotes a prediction of the HT solution at $t^{(k)}$ obtained using the first order scheme, included so that the augmented bases and core tensors contain information over the entire subinterval $[t^n,t^{(k)}]$. The computational cost of the reduced augmentation procedure in equation \eqref{eq:redaug_HT} is negligible relative to the entire algorithm \cite{KressnerTobler2014} which is dominated by the Sylvester steps. Defining $\mathbf{K}_1^{(k)}$ as in equation \eqref{eq:K1_K2} but at $t^{(k)}$ and onto the approximate bases, we can rearrange and project equation \eqref{eq:impuprk} similarly to equation \eqref{eq:impup} in Subsection \ref{subsec:firstorder} to obtain the Sylvester equation $\mathbf{A}\mathbf{K}_1^{(k)}+\mathbf{K}_1^{(k)}\mathbf{B}=\mathbf{C}$, where $\mathbf{A}=\mathbf{I}_{N_1}-a_{kk}\Delta t\mathbf{F}_1$,

\begin{align}
\begin{split}
    \mathbf{B} =& - a_{kk}\Delta t\Big( (\mathbf{B}_{34}^{\star,(k)})^T((\mathbf{F}_4\mathbf{V}_4^{\star,(k)})^T\mathbf{V}_4^{\star,(k)}\otimes\mathbf{I}_{r_3})\mathbf{B}_{34}^{\star,(k)}\otimes\mathbf{I}_{r_2}\\
    &+ (\mathbf{B}_{34}^{\star,(k)})^T(\mathbf{I}_{r_4}\otimes(\mathbf{F}_3\mathbf{V}_3^{\star,(k)})^T\mathbf{V}_3^{\star,(k)})\mathbf{B}_{34}^{\star,(k)}\otimes\mathbf{I}_{r_2}\\
    &+(\mathbf{I}_{r_{34}}\otimes (\mathbf{F}_2\mathbf{V}_2^{\star,(k)})^T\mathbf{V}_2^{\star,(k)}
     \Big),
\end{split}
\end{align}

\begin{align}\label{eq:impup_mat_rightmultrk}
\begin{split}
\mathbf{C}
&=
\mathbf{V}_1^{n}
\tilde{\mathbf{B}}_{12,(a)}^{n}
\big(
\mathbf{B}_{34}^{n,T}(\mathbf{V}_4^{n,T}\mathbf{V}_{4}^{\star,(k)} \otimes \mathbf{V}_3^{n,T}\mathbf{V}_{3}^{\star,(k)}) \mathbf{B}_{34}^{\star,(k)} \otimes \mathbf{V}_2^{n,T}\mathbf{V}_2^{\star,(k)}
\big)\\
&+ \Delta t \sum_{\ell=1}^{k-1}a_{k\ell}\Big\{
 \mathbf{F}_1 \mathbf{V}_1^{(\ell)}
\tilde{\mathbf{B}}_{12,(a)}^{(\ell)}
\big(
\mathbf{B}_{34}^{(\ell),T}(\mathbf{V}_4^{(\ell),T}\mathbf{V}_{4}^{\star,(k)} \otimes \mathbf{V}_3^{(\ell),T}\mathbf{V}_{3}^{\star,(k)}) \mathbf{B}_{34}^{\star,(k)} \otimes \mathbf{V}_2^{(\ell),T}\mathbf{V}_2^{\star,(k)}
\big)\\
&\qquad\quad\qquad+\mathbf{V}_1^{(\ell)}
\tilde{\mathbf{B}}_{12,(a)}^{(\ell)}
\big(
\mathbf{B}_{34}^{(\ell),T}(\mathbf{V}_4^{(\ell),T}\mathbf{V}_{4}^{\star,(k)} \otimes \mathbf{V}_3^{(\ell),T}\mathbf{V}_{3}^{\star,(k)}) \mathbf{B}_{34}^{\star,(k)} \otimes (\mathbf{F}_2\mathbf{V}_2^{(\ell)})^T\mathbf{V}_2^{\star,(k)}
\big)\\
&\qquad\quad\qquad+ \mathbf{V}_1^{(\ell)}
\tilde{\mathbf{B}}_{12,(a)}^{(\ell)}
\big(\mathbf{B}_{34}^{(\ell),T}
(\mathbf{V}_4^{(\ell),T}\mathbf{V}_4^{\star,(k)}
\otimes
(\mathbf{F}_3\mathbf{V}_3^{(\ell)})^{T}\mathbf{V}_3^{\star,(k)}
\big)
\mathbf{B}_{34}^{\star,(k)}
\otimes
\mathbf{V}_2^{(\ell),T}\mathbf{V}_2^{\star,(k)}\big)\\
&\qquad\quad\qquad+ \mathbf{V}_1^{(\ell)}
\tilde{\mathbf{B}}_{12,(a)}^{(\ell)}
\big(\mathbf{B}_{34}^{(\ell),T}
(\mathbf{F}_4\mathbf{V}_4^{(\ell)})^T\mathbf{V}_4^{\star,(k)}
\otimes
\mathbf{V}_3^{(\ell),T}\mathbf{V}_3^{\star,(k)}
\big)
\mathbf{B}_{34}^{\star,(k)}
\otimes
\mathbf{V}_2^{(\ell),T}\mathbf{V}_2^{\star,(k)}\big)\Big\}.
\end{split}
\end{align}

The K-step updated basis $\mathbf{V}_1^{\ddagger,(k)}$ is obtained by computing the orthonormal portion of the reduced QR factorization of $\mathbf{K}_1^{(k)}$ as in Subsection $\ref{subsec:firstorder}$. We now reduce the augmented matrix 
\begin{equation}
    [\mathbf{V}_1^{\ddagger,(k)},\mathbf{V}_1^{(k-1)},...,\mathbf{V}_1^{(1)},\mathbf{V}_1^{(n)}]
\end{equation}
and to obtain $\hat{\mathbf{V}}_1^{(k)}$. The other 1D bases $\hat{\mathbf{V}}_i^{(k)}$ ($i=2,3,4$) are obtained similarly. For the B12 step, we perform a Galerkin projection onto the updated bases and let
\begin{equation}
\operatorname{vec}(\hat{\tilde{\mathcal B}}^{(k)}_{12} )\ \coloneq \left((\mathbf{V}_{34}^{\star,(k)})^{T} \mathbf{V}_{34}^{{(k)}}
    \otimes (\hat{\mathbf{V}}_{2}^{(k)})^T \mathbf{V}_{2}^{{(k)}}
    \otimes (\hat{\mathbf{V}}_{1}^{(k)})^T\mathbf{V}_{1}^{{(k)}}
  \right)
\operatorname{vec}(\tilde{\mathcal{B}}_{12}^{{(k)}}). 
\end{equation}

We then vectorize equation \eqref{eq:impuprk}, and multiply on the left by $(\mathbf{V}_{34}^{\star,(k)}\otimes\hat{\mathbf{V}}_2^{(k)}\otimes\hat{\mathbf{V}}_1^{(k)})^T$ to get the third order tensor linear equation
\begin{equation} \label{eq:B12tilhatrk}
\begin{aligned}
\left\{
(\mathbf{B}_{34}^{\star,(k)})^T
\Big(
  -a_{kk} \Delta t (\hat{\mathbf{V}}_4^{(k)})^T
  (\mathbf{F}_4 \hat{\mathbf{V}}_4^{(k)})
  \otimes \mathbf{I}_{\hat{r}_3^{(k)}}
\Big)
\mathbf{B}_{34}^{\star,(k)} \otimes \mathbf{I}_{\hat{r}_2^{(k)}} \otimes \mathbf{I}_{\hat{r}_1^{(k)}}
\right. \\ \left.
+ (\mathbf{B}_{34}^{\star,(k)})^T
\Big(
  \mathbf{I}_{\hat{r}_4^{(k)}}
  \otimes
  -a_{kk} \Delta t (\hat{\mathbf{V}}_4^{(k)})^T
  (\mathbf{F}_4 \hat{\mathbf{V}}_4^{(k)})
\Big)
\mathbf{B}_{34}^{\star,(k)}
\right) \otimes \mathbf{I}_{\hat{r}_2^{(k)}} \otimes \mathbf{I}_{\hat{r}_1^{(k)}}\\ \left.
+ \mathbf{I}_{\hat{r}_{34}^{(k)}}
\otimes
\Big(
  -a_{kk} \Delta t (\hat{\mathbf{V}}_2^{(k)})^T
  (\mathbf{F}_2 \hat{\mathbf{V}}_2^{(k)})
\Big)
\otimes \mathbf{I}_{\hat{r}_1^{(k)}}
\right. \\ \left.
+ \mathbf{I}_{\hat{r}_{34}^{(k)}}
\otimes \mathbf{I}_{\hat{r}_2^{(k)}}
\otimes
\Big(
  -a_{kk} \Delta t (\hat{\mathbf{V}}_1^{(k)})^T
  (\mathbf{F}_1 \hat{\mathbf{V}}_1^{(k)})
\Big)
\right\}
\operatorname{vec}(\hat{\tilde{\mathcal{B}}}_{12}^{(k)})
\\= \operatorname{vec}(\mathcal{S}_{12}^{(k-1)}),
\end{aligned}
\end{equation}
where  
\begin{equation}
\label{eq:b12_rk_rhs}
\begin{aligned}
&\quad\mathcal{S}_{12}^{(k-1)} 
=\tilde{\mathcal{B}}_{12}^{n}
\times_a \hat{\mathbf V}_1^{(k),T}\mathbf{V}_1^{n}
\times_b \hat{\mathbf V}_2^{(k),T}\mathbf{V}_2^{n}
\times_c \hat{\mathbf B}_{34}^{(k),T}
\big(
\hat{\mathbf V}_{4}^{(k),T}\mathbf{V}_{4}^{n}
\otimes
\hat{\mathbf V}_{3}^{(k),T}\mathbf{V}_{3}^{n}
\big)\mathbf B_{34}^{n}
\\
&+\Delta t \sum_{\ell=1}^{k-1}a_{k\ell}\Big\{
\tilde{\mathcal{B}}_{12}^{(\ell)}
\times_a \hat{\mathbf V}_1^{(k),T}(\mathbf F_1\mathbf{V}_1^{(\ell)})
\times_b \hat{\mathbf V}_2^{(k),T}\mathbf{V}_2^{(\ell)}
\times_c \hat{\mathbf B}_{34}^{(k),T}
\big(
\hat{\mathbf V}_{4}^{(k),T}\mathbf{V}_{4}^{(\ell)}
\otimes
\hat{\mathbf V}_{3}^{(k),T}\mathbf{V}_{3}^{(\ell)}
\big)\mathbf B_{34}^{(\ell)}
\\
&\qquad\qquad+
\tilde{\mathcal{B}}_{12}^{(\ell)}
\times_a \hat{\mathbf V}_1^{(k),T}\mathbf{V}_1^{(\ell)}
\times_b \hat{\mathbf V}_2^{(k),T}(\mathbf F_2\mathbf{V}_2^{(\ell)})
\times_c \hat{\mathbf B}_{34}^{(k),T}
\big(
\hat{\mathbf V}_{4}^{(k),T}\mathbf{V}_{4}^{(\ell)}
\otimes
\hat{\mathbf V}_{3}^{(k),T}\mathbf{V}_{3}^{(\ell)}
\big)\mathbf B_{34}^{(\ell)}
\\
&\qquad\qquad+
\tilde{\mathcal{B}}_{12}^{(\ell)}
\times_a \hat{\mathbf V}_1^{(k),T}\mathbf{V}_1^{(\ell)}
\times_b \hat{\mathbf V}_2^{(k),T}\mathbf{V}_2^{(\ell)}
\times_c \hat{\mathbf B}_{34}^{(k),T}
\big(
\hat{\mathbf V}_{4}^{(k),T}\mathbf{V}_{4}^{(\ell)}
\otimes
\hat{\mathbf V}_{3}^{(k),T}(\mathbf F_3\mathbf{V}_3^{(\ell)})
\big)\mathbf B_{34}^{(\ell)}
\\
&\qquad\qquad+
\tilde{\mathcal{B}}_{12}^{(\ell)}
\times_a \hat{\mathbf V}_1^{(k),T}\mathbf{V}_1^{(\ell)}
\times_b \hat{\mathbf V}_2^{(k),T}\mathbf{V}_2^{(\ell)}
\times_c \hat{\mathbf B}_{34}^{(k),T}
\big(
\hat{\mathbf V}_{4}^{(k),T}(\mathbf F_4\mathbf{V}_4^{(\ell)})
\otimes
\hat{\mathbf V}_{3}^{(k),T}\mathbf{V}_{3}^{(\ell)}
\big)\mathbf B_{34}^{(\ell)}
\Big\}.
\end{aligned}
\end{equation}

We solve equation \eqref{eq:B12tilhatrk} similarly to equation \eqref{eq:impup_vec_proj_lin} using the direct solver of Simoncini \cite{Simoncini2020}. Updating $\mathcal{B}_{34}$ follows similarly. Updating the core tensor at the leaf node $\mathcal{B}_{1234}$ follows a similar procedure to the backward Euler case in equation \eqref{eq:b1234} with righthand side terms similar to equation \eqref{eq:b12_rk_rhs}. We finally truncate using with a tolerance of $\varepsilon=10^{-6}$ to obtain our solution in HT format.

\begin{rem}\label{rem:highD}
Higher-dimensional extension follows the same procedure, evolving from the leaf bases up to the root node. For example, computing the solution to a 6D equation with solutions decomposed according to the tree shown in Figure \ref{fig:6Dtree} involves a set of K-steps for the leaf bases on nodes $3,4,5,6$ followed by an evolution of the core tensors at the $34$ and $56$ nodes. A set of K-steps at nodes $1,2$ follows, then the cores at the $12$ and $3456$ nodes are updated, and we finally proceed to evolve the core tensor at the root node.
\begin{figure}[H]
    \centering

\begin{forest}
  [\(\mathcal{B}_{123456}\)
    [\(\mathcal{B}_{12}\)
      [\(\mathbf{V}_1\)]
      [\(\mathbf{V}_2\)]
    ]
    [\(\mathcal{B}_{3456}\)
      [\(\mathcal{B}_{34}\)
        [\(\mathbf{V}_3\)]
        [\(\mathbf{V}_4\)]
      ]
      [\(\mathcal{B}_{56}\)
        [\(\mathbf{V}_5\)]
        [\(\mathbf{V}_6\)]
      ]
    ]
  ]
\end{forest}
   \caption{A hierarchical Tucker decomposition of a 6D tensor.}
    \label{fig:6Dtree}
\end{figure}
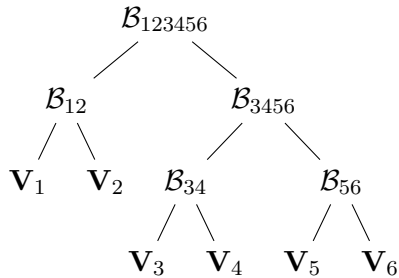
\end{rem}

\section{Numerical Experiments}
\label{sec:numexp}
We now numerically validate the high-order accuracy (in $L^1$) of the proposed algorithm to solve the diffusion equation \eqref{eq:4dheat}. Rank plots illustrate the method’s ability to capture low-rank solution structures through a collection of benchmark diffusion coefficients. In each example, we assume a uniform mesh with ($N_{x_i}=60$) grid points for $i = {1,2,3,4}$. On this uniform mesh, one-dimensional spatial derivatives are discretized using spectral methods, with differentiation matrices given in \cite{Trefethan}. Consequently, the temporal error is the only source of error. Periodic boundary conditions are assumed throughout, but the framework can be extended to accommodate alternative boundary conditions by employing the corresponding differentiation matrices. For problems involving solutions that decay at infinity, such as Gaussian distribution functions, the computational domain is chosen sufficiently large so that the solution remains smooth near the boundaries, thereby justifying the use of spectral methods with periodic boundary conditions. 
Since we only consider diffusion equations we define the time step size $\Delta t=\lambda\Delta x/4$, where $\Delta x=\Delta x_1=\Delta x_2=\Delta x_3=\Delta x_4$.

\subsection{Accuracy Test}
For the accuracy test, we compute the $L^1$ error for $\lambda \in \{0.1, 1.05,...,2.95,3\}$, $d_1 = d_2 = d_3 =  d_4 = 0.1$, final time $T_f = 0.5$, and domain $\Omega=[0,2\pi]^4$. 
We divide the $L^1$ error by the measure of the domain $|\mathbf{\Omega}|$ so that it is comparable to the $L^{\infty}$ norm. We let the initial condition be the first three Fourier modes,
\begin{equation}\label{eq:accic}
    u_0({x}) = \sum_{k=1}^{3} 
\prod_{i=1}^{4} \sin(k x_i).
\end{equation}

We observe the expected temporal accuracy and remain rank-3 in each dimension using backward Euler, DIRK2, and DIRK3. By equation \eqref{eq:HT4complx} our rank-3 problem observes a storage savings of approximately 17,400x compared to naively storing the full $N^4$ solution.

\begin{minipage}{0.50\textwidth}
\begin{figure}[H]
    \centering
    \includegraphics[width=\linewidth]{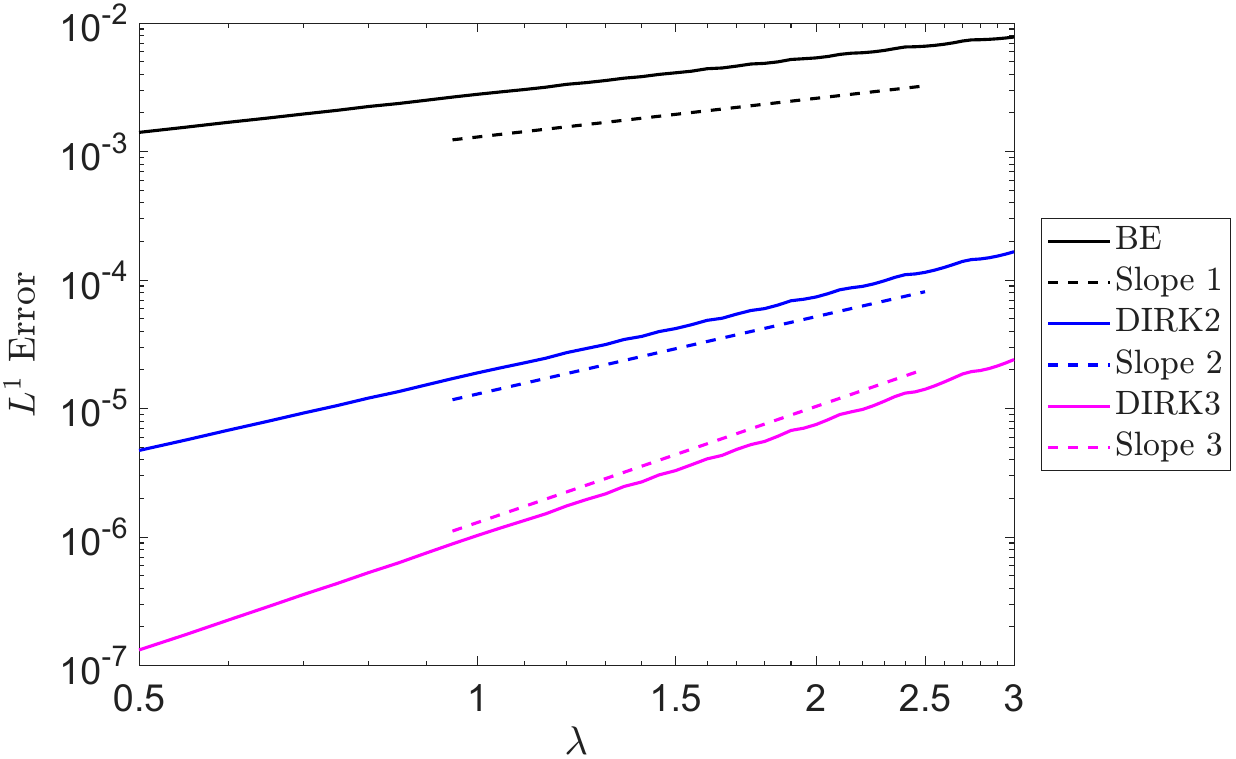}
    \caption{Error plot for equation \eqref{eq:4dheat} using initial condition \eqref{eq:accic} using Backward Euler, DIRK2, and DIRK3}
    \label{fig:accplot}
\end{figure}
\end{minipage}
\begin{minipage}{0.05\textwidth}
\end{minipage}
\begin{minipage}{0.41\textwidth}
    \begin{figure}[H]
        \centering
        \includegraphics[width=\linewidth]{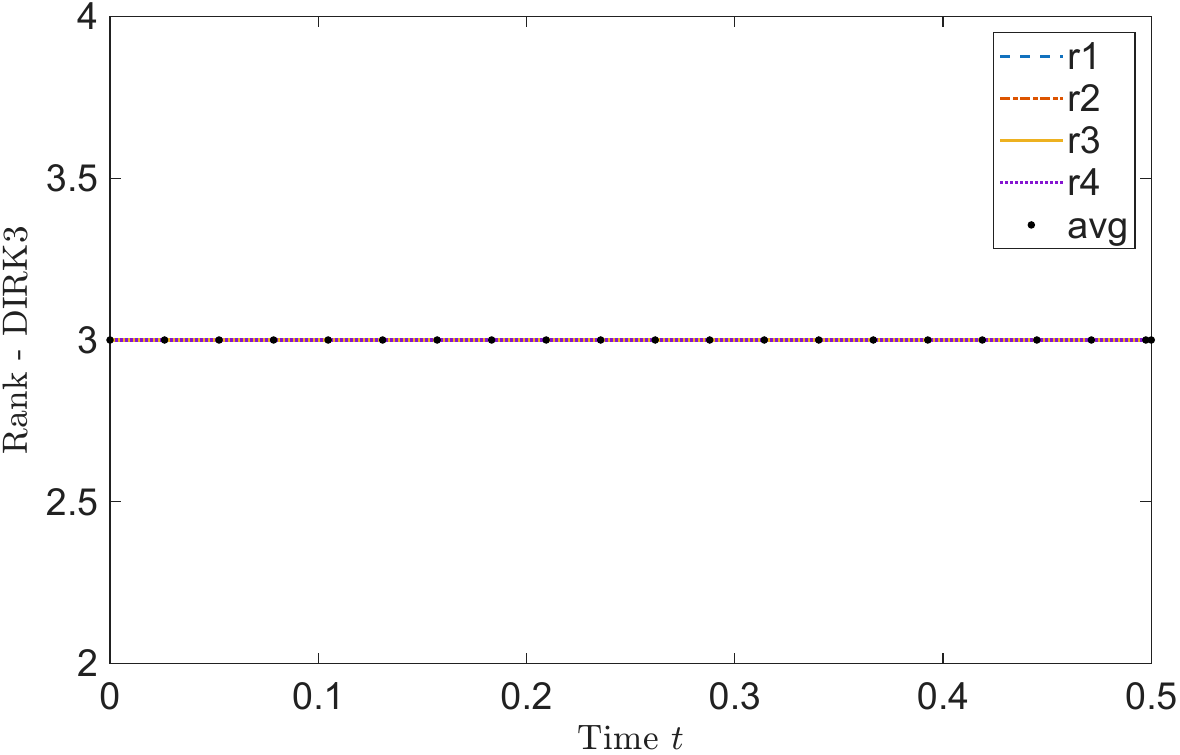}
        \ \\
        \ \\
        \caption{Rank Plot for equation \eqref{eq:4dheat} using initial condition \eqref{eq:accic} using DIRK3}
        \label{fig:rankplotacc}
    \end{figure}
\end{minipage}

\subsection{Rank Tests}
The following tests show how well the algorithm captures the physical rank of the solution. We consider the rank-3 initial condition 
\begin{equation}\label{rank3ic}
u(\mathbf{x},t=0)=0.8\exp\left(-15|\mathbf{x}-6.5|^2\right)+0.5\exp\left(-15|\mathbf{x}-7.5|^2\right)+1.2\exp(-15|\mathbf{x}-4.5|^2).
\end{equation}

The numerical domain $\Omega=[0,14]^4$ is large enough so that spectral methods can be used while assuming periodic boundary conditions for the initial condition composed of Gaussian functions. We set $\lambda = 1, T_f = 15$.
We test different sets of diffusion coefficients in equation \eqref{eq:4dheat}. Due to rapid changes in the solution at short times from diffusion, we expect the rank to initially increase. However as the solution evolves, we expect the ranks to decrease according to the diffusion coefficients.

\subsubsection*{Constant Coefficients}
\begin{figure}[H]
    \centering
    \begin{subfigure}[b]{0.31\textwidth}
        \centering
        \includegraphics[width=\textwidth]{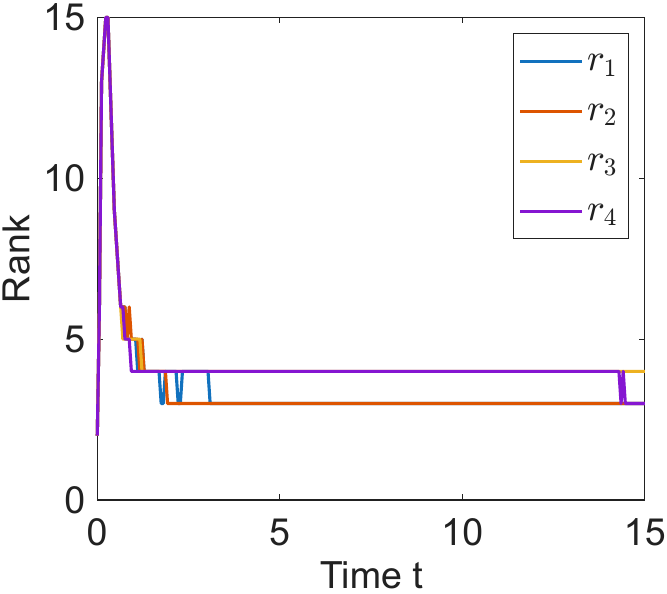}
        \caption{Backward Euler}
        \label{fig:beluer_const_leaf}
    \end{subfigure}
    \hfill
    \begin{subfigure}[b]{0.31\textwidth}
        \centering
        \includegraphics[width=\textwidth]{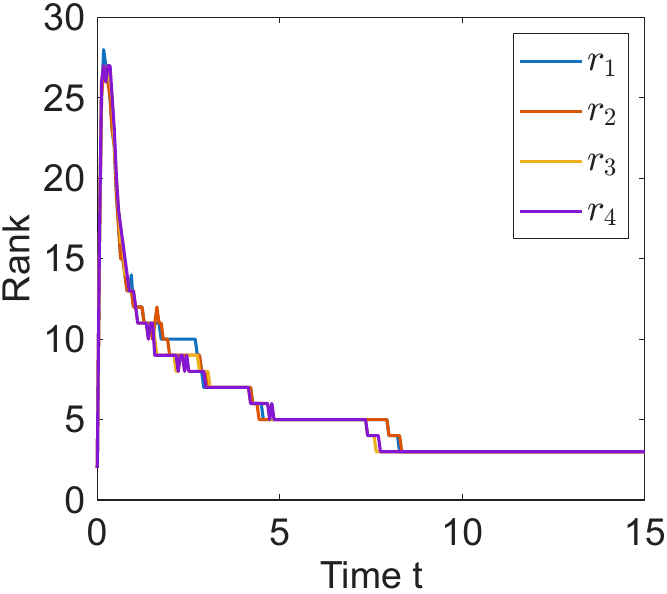}
        \caption{DIRK 2}
        \label{fig:dirk2_const_leaf}
    \end{subfigure}
    \hfill
    \begin{subfigure}[b]{0.31\textwidth}
        \centering
        \includegraphics[width=\textwidth]{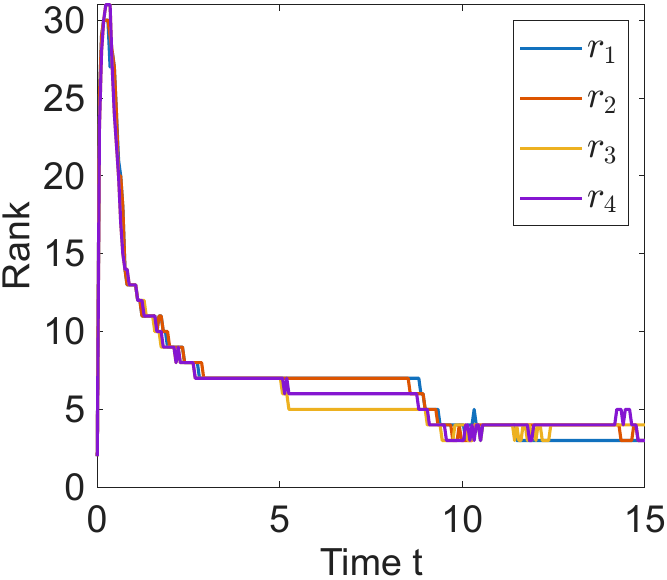}
        \caption{DIRK 3}
        \label{fig:dirk3_const_leaf}
    \end{subfigure}
    \caption{Rank plots for the leafs, $r_1,r_2,r_3,r_4$ for constant diffusion.}
    \label{fig:solutions_all_leaf}
\end{figure}

\begin{figure}[H]
    \centering
    \begin{subfigure}[b]{0.31\textwidth}
        \centering
        \includegraphics[width=\textwidth]{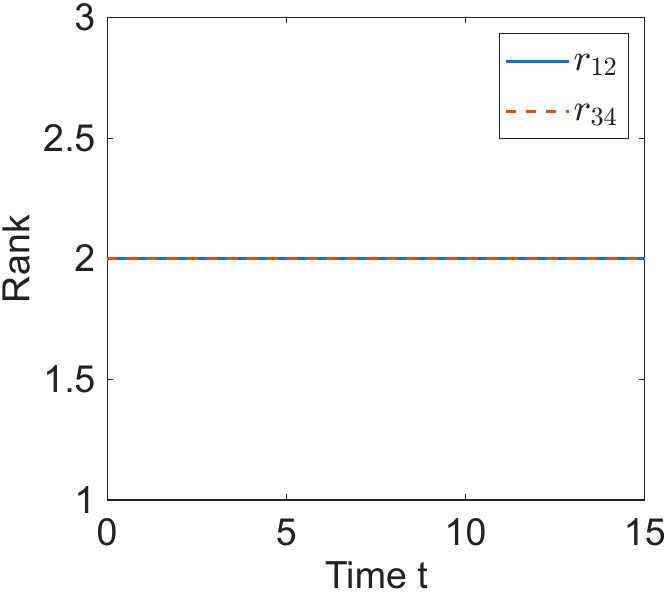}
        \caption{Backward Euler}
        \label{fig:beuler_const_core}
    \end{subfigure}
    \hfill
    \begin{subfigure}[b]{0.31\textwidth}
        \centering
        \includegraphics[width=\textwidth]{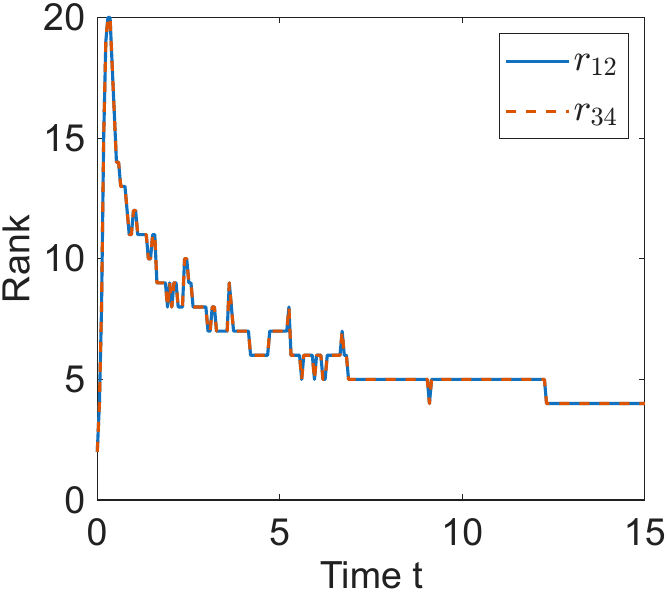}
        \caption{DIRK 2}
        \label{fig:dirk2_const_core}
    \end{subfigure}
    \hfill
    \begin{subfigure}[b]{0.31\textwidth}
        \centering
        \includegraphics[width=\textwidth]{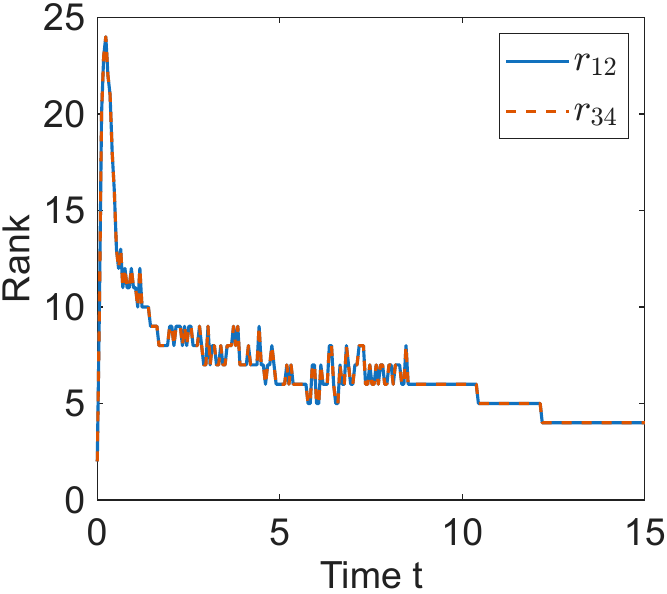}
        \caption{DIRK 3}
        \label{fig:dirk3_const_core}
    \end{subfigure}
    \caption{Rank plots for the $r_{12},r_{34}$ for constant diffusion.}
    \label{fig:const_all_core}
\end{figure}

$d_1 = d_2 = d_3 = d_4 = 1$.
For the case of constant diffusion coefficients, we expect a sharp initial rise in rank due to short-time diffusive dynamics, followed by a gradual decrease. We observe this behavior in each of the \textit{leaf modes} ($r_1,r_2,r_3,r_4$) for backward Euler, DIRK2, and DIRK3, as seen in Figure \ref{fig:solutions_all_leaf}. However, for the backward Euler test, $r_{12}$ and $r_{34} $ remain rank-2. Crucially, the backward Euler test only uses information at time $t^n$ to predict bases that are changing very quickly at short times. Thus the evolution equations for the core tensors do not see any information from intermediate timesteps, and so over-truncates the solution. This lack of information explains why $r_1,r_2,r_3,r_4$ in the backward Euler scheme still increase initially, but sharply decrease faster than its higher order counterparts. But, this idea is best seen in the intermediate ranks $r_{12},r_{34}$ in which the rank is constantly 2, which is not physically expected, as seen in Figure \ref{fig:const_all_core}.


\subsubsection*{Sinusoidal Coefficients}
\begin{figure}[H]
    \centering
    \begin{subfigure}[b]{0.31\textwidth}
        \centering
        \includegraphics[width=\textwidth]{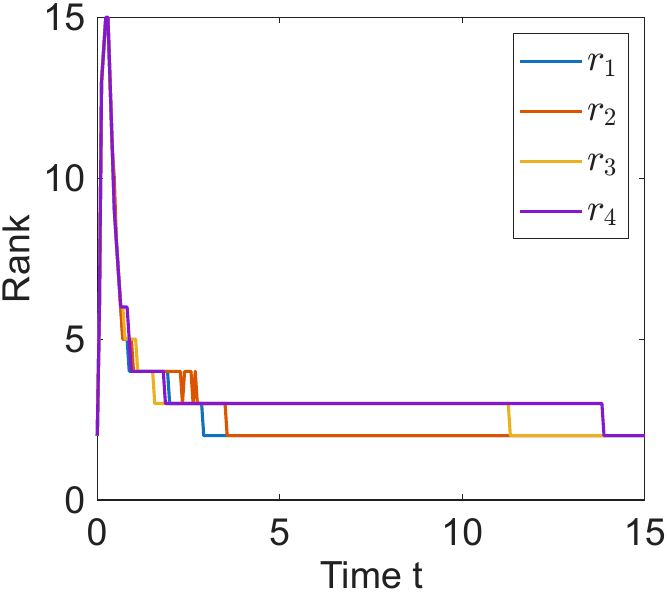}
        \caption{Backward Euler}
        \label{fig:beuler_square_leaf}
    \end{subfigure}
    \hfill
    \begin{subfigure}[b]{0.31\textwidth}
        \centering
        \includegraphics[width=\textwidth]{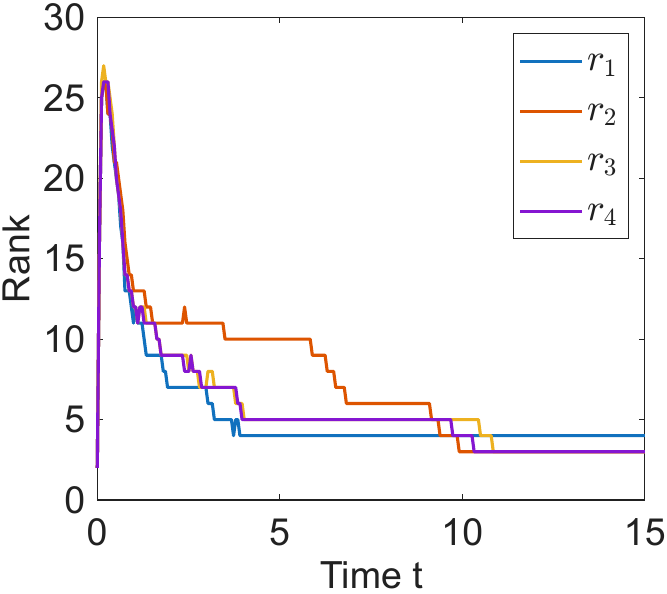}
        \caption{DIRK 2}
        \label{fig:dirk2_square_leaf}
    \end{subfigure}
    \hfill
    \begin{subfigure}[b]{0.31\textwidth}
        \centering
        \includegraphics[width=\textwidth]{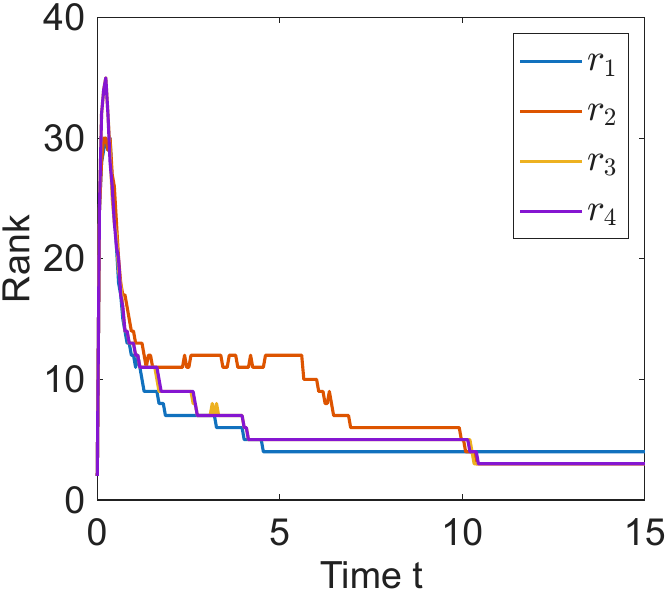}
        \caption{DIRK 3}
        \label{fig:dirk3_square_leaf}
    \end{subfigure}
    \caption{Rank plots for the leafs, $r_1,r_2,r_3,r_4$ for sinusoidal diffusion.}
    \label{fig:sin_all_leaf}
\end{figure}

$d_1 = 1+\sin(2\pi t/T_f), d_2 = 1-\sin(2\pi t/T_f) , d_3 = d_4 = 1$. Setting $T_f=15$, we see that for DIRK2 and DIRK3, $r_1$ diffuses more rapidly than all other ranks and $r_2$ more slowly until $T_f/2=7.5$, as expected from the diffusion coefficients. And as $d_2$ decreases then increases over $0\leq t\leq 7.5$, the rate at which $r_2$ decreases is as expected, with $dr_2/dt=0$ around $t=3.75$ since $d_2(3.75)=0$. However, this rank behavior is not detectable in the backward Euler test, as seen in Figure \ref{fig:sin_all_leaf}. For the same underlying reason as the backward Euler case with constant coefficients, the rank of the solution at nodes $12$ and $34$ do not exhibit the physically expected behavior. This underscores the importance of the higher order methods studied for rank capture in the solution. 

\begin{figure}[H]
    \centering
    \begin{subfigure}[b]{0.31\textwidth}
        \centering
        \includegraphics[width=\textwidth]{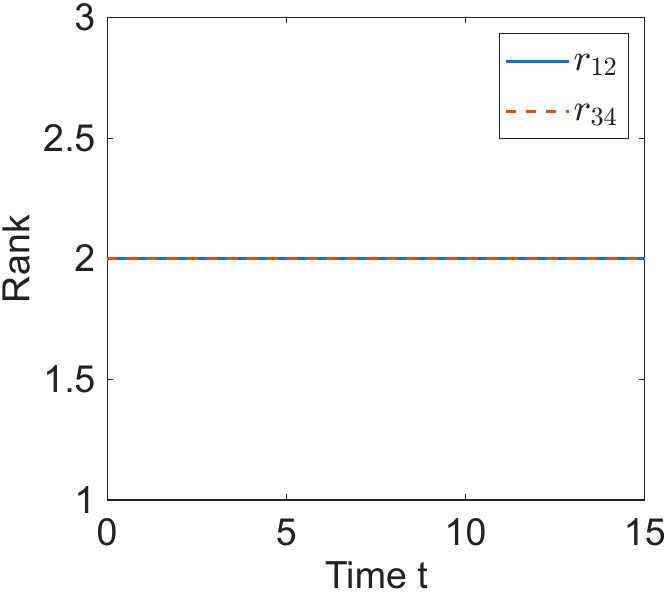}
        \caption{Backward Euler}
        \label{fig:sin_beuler_core}
    \end{subfigure}
    \hfill
    \begin{subfigure}[b]{0.31\textwidth}
        \centering
        \includegraphics[width=\textwidth]{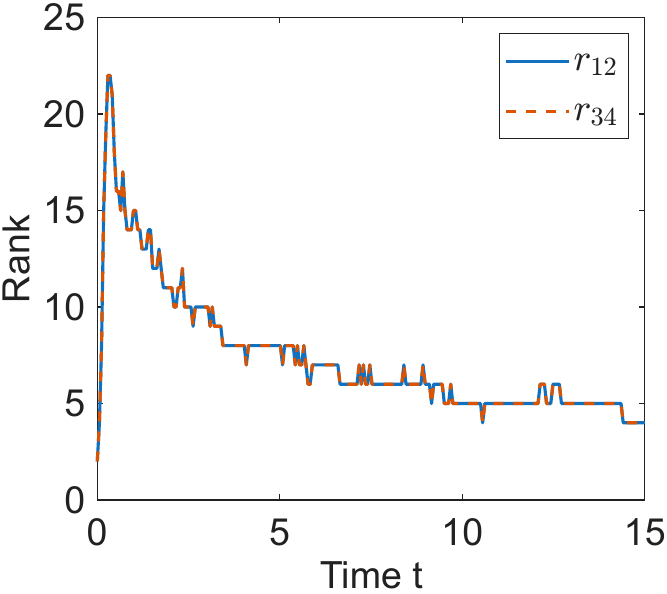}
        \caption{DIRK 2}
        \label{fig:sin_dirk2_core}
    \end{subfigure}
    \hfill
    \begin{subfigure}[b]{0.31\textwidth}
        \centering
        \includegraphics[width=\textwidth]{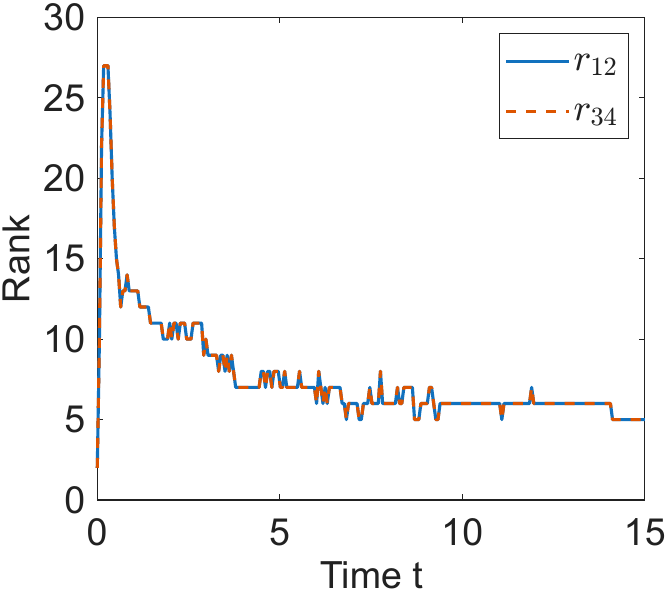}
        \caption{DIRK 3}
        \label{fig:sin_dirk3_core}
    \end{subfigure}
    \caption{Rank plots for the $r_{12},r_{34}$ for sinusoidal diffusion.}
    \label{fig:sin_all_core}
\end{figure}

\subsubsection*{Square Wave Coefficients}
\begin{figure}[!htbp]
    \centering
    \begin{subfigure}[b]{0.31\textwidth}
        \centering
        \includegraphics[width=\textwidth]{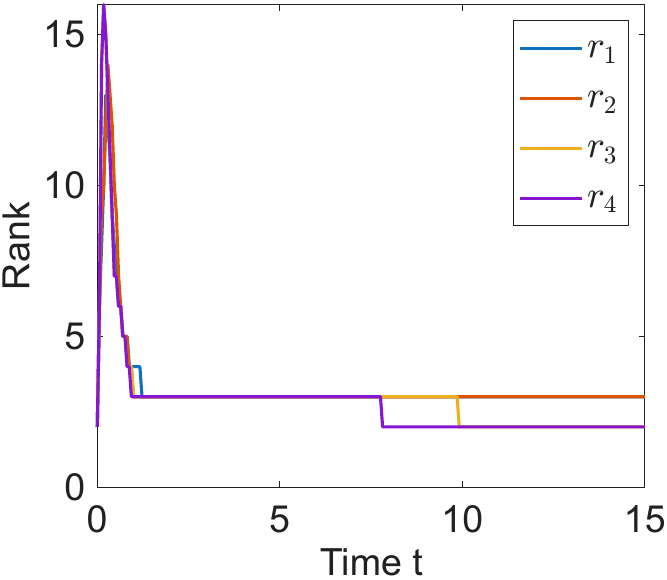}
        \caption{Backward Euler}
        \label{fig:square_beuler_leaf}
    \end{subfigure}
    \hfill
    \begin{subfigure}[b]{0.31\textwidth}
        \centering
        \includegraphics[width=\textwidth]{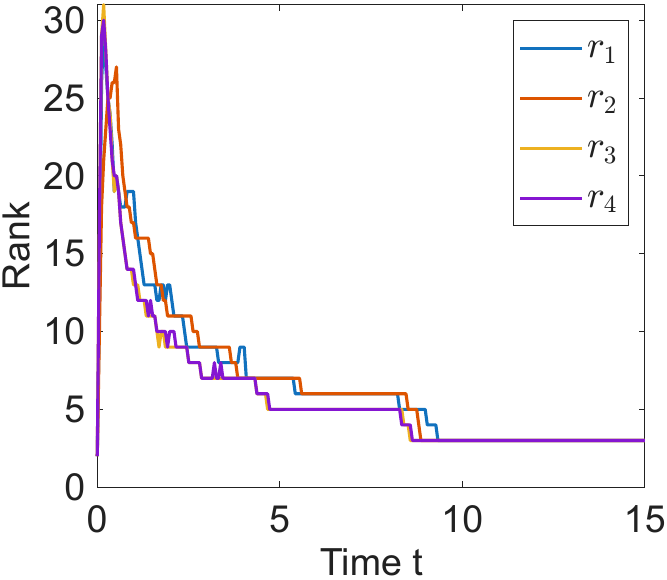}
        \caption{DIRK 2}
        \label{fig:square_dirk2_leaf}
    \end{subfigure}
    \hfill
    \begin{subfigure}[b]{0.31\textwidth}
        \centering
        \includegraphics[width=\textwidth]{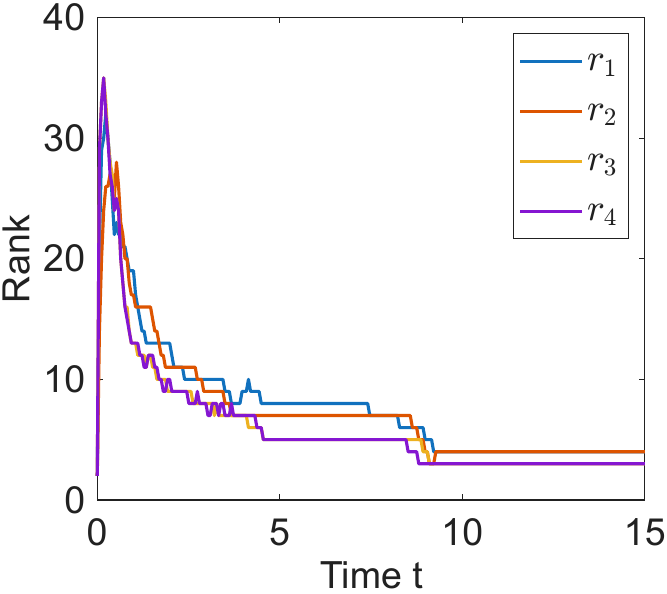}
        \caption{DIRK 3}
        \label{fig:square_dirk3_leaf}
    \end{subfigure}
    \caption{Rank plots for the leafs, $r_1,r_2,r_3,r_4$ for square wave diffusion.}
    \label{fig:square_all_leaf}
\end{figure}

$d_1 = 1.05+\operatorname{sign}(\sin(2\pi t))/2, d_2 = 1.05-\operatorname{sign}(\sin(2\pi t))/2, d_3 = 1, d_4 = 1$.
We see in Figure \ref{fig:square_all_leaf} that $r_1$ and $r_2$ decrease more rapidly when their respective diffusion coefficient is higher, as expected. We also see that $r_2$ reaches a lower maximum rank than $r_1$, because its diffusion coefficient begins lower, as expected. We continue to see similar nonphysical behavior in the backward Euler test as before. 
\begin{figure}[H]
    \centering
    \begin{subfigure}[b]{0.31\textwidth}
        \centering
        \includegraphics[width=\textwidth]{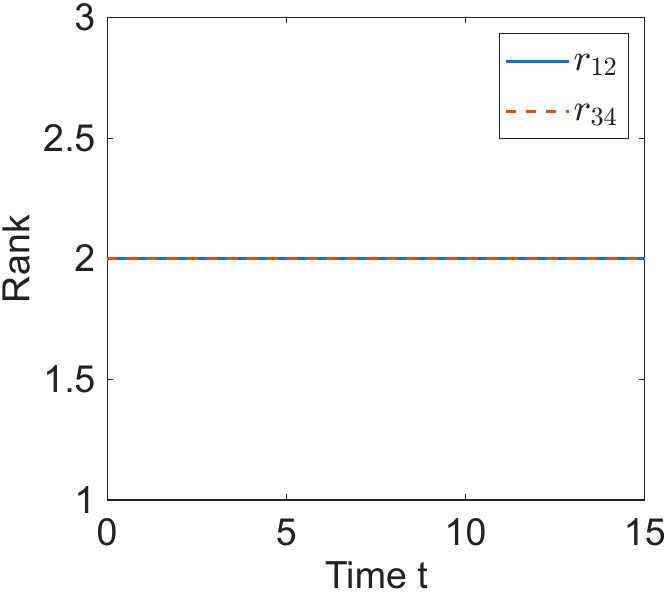}
        \caption{Backward Euler}
        \label{fig:square_beuler_core}
    \end{subfigure}
    \hfill
    \begin{subfigure}[b]{0.31\textwidth}
        \centering
        \includegraphics[width=\textwidth]{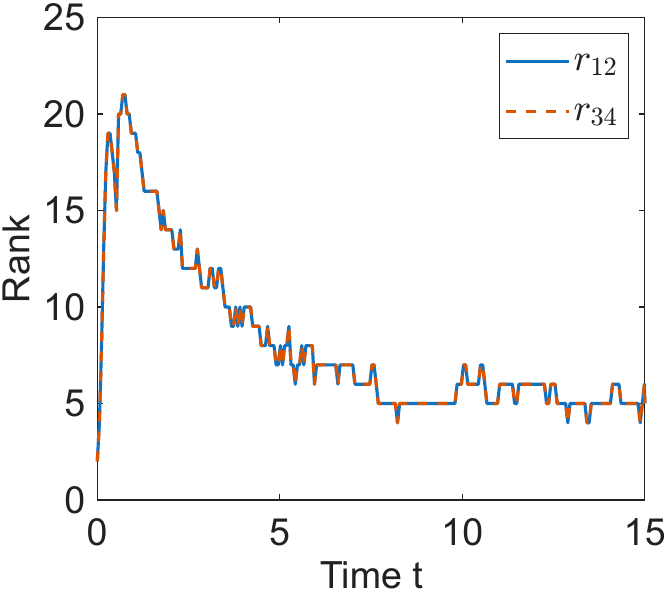}
        \caption{DIRK 2}
        \label{fig:square_dirk2_core}
    \end{subfigure}
    \hfill
    \begin{subfigure}[b]{0.31\textwidth}
        \centering
        \includegraphics[width=\textwidth]{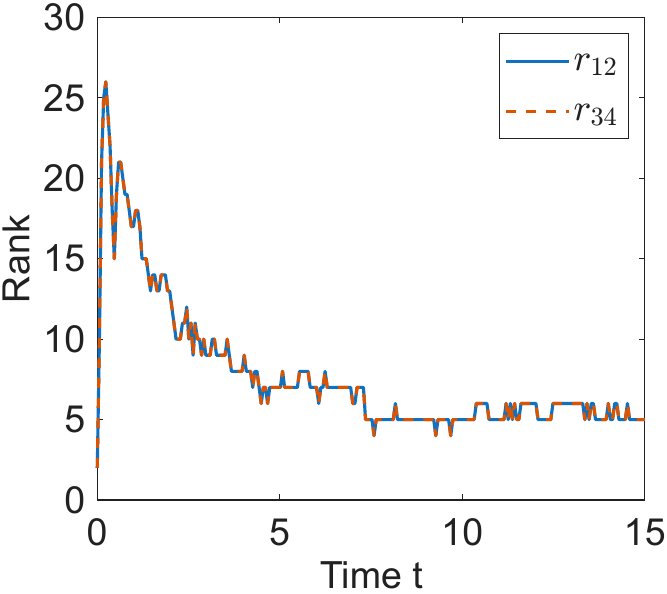}
        \caption{DIRK 3}
        \label{fig:square_dirk3_core}
    \end{subfigure}
    \caption{Rank plots for the $r_{12},r_{34}$ for square wave diffusion.}
    \label{fig:square_all_core}
\end{figure}

\section{Conclusion}
\label{sec:conclusions}
In this paper, we present a reduced augmentation implicit low-rank (RAIL) algorithm for solving high dimensional diffusion equations. The tensor equations were projected onto an enriched subspace formed by combining a low-order prediction (or updated) basis with bases from earlier RK stages through a reduced augmentation strategy, thereby permitting the use of high-order implicit schemes. Diffusion equations were discretized in space using the hierarchical Tucker decomposition, and spectral discretizations for spatial derivatives. Several results support its efficacy at capturing low-rank structure in high dimensional diffusion equations. We presented results for the $d=4$ case, but this algorithm can be extended easily higher dimensional problems, and its proof of concept sets the stage for a series of algorithms that efficiently solve other high dimensional PDEs. Current work focuses on extending this algorithm to advection-diffusion equations using implicit-explicit methods.

\subsection{Data Availability} The source code generated and analyzed in the current study is available upon reasonable request.

\subsection{Competing Interests} The author has no competing interests.

\subsection{Acknowledgments} The author is grateful for the support and mentorship of Dr. Joseph Nakao. The author would also like to thank Dr. Gianluca Ceruti for helpful input and discussions.

\appendix
\section{Diagonally Implicit Runge-Kutta (DIRK) Butcher Tableaus \cite{ascher1997implicit}}
\begin{table}[H]
\centering
\begin{subtable}{0.45\textwidth}
\centering
\begin{equation*}
\begin{aligned}
\begin{array}{c|cc}
\nu & \nu & 0 \\
1   & 1-\nu & \nu \\
\hline
    & 1-\nu & \nu
\end{array}
\\
\nu = 1-\sqrt{2}/2
\end{aligned}
\end{equation*}

\end{subtable}
\hfill
\begin{subtable}{0.45\textwidth}
\centering
\begin{equation*}
\begin{aligned}
\begin{array}{c|ccc}
\nu & \nu & 0 & 0 \\
\frac{1+\nu}{2} & \frac{1-\nu}{2} & \nu & 0 \\
1 & \beta_1 & \beta_2 & \nu \\
\hline
  & \beta_1 & \beta_2 & \nu
\end{array}
\\
\nu \approx 0.435866521508459\\
\beta_1 = -(3/2)\nu^2 + 4\nu - 1/4\\
\beta_2 = (3/2)\nu^2 - 5\nu + 5/4
\end{aligned}
\end{equation*}
\end{subtable}
\caption{Butcher tables for stiffly accurate DIRK methods. (Left) Second-order DIRK method. (Right) Third-order DIRK method.}
\label{tab:DIRK_tables}
\end{table}



\medskip

\printbibliography

\end{document}